\documentclass[11pt]{article}

\usepackage{latexsym}
\usepackage{amssymb}
\usepackage{amsthm}
\usepackage{amscd}
\usepackage{amsmath}
\usepackage[all]{xy}
\usepackage{graphicx}

\newtheorem{theorem}{Theorem}[section]

\newtheorem{lemma}{Lemma}[section]

\newtheorem{proposition}{Proposition}[section]
\newtheorem{corollary}{Corollary}[section]

\xyoption{dvips}

\numberwithin{equation}{section}

\newcommand{\FF}{\mathbb{F}}    \newcommand{\CC}{\mathbb{C}}    
 \def\dim{\mathrm{dim}}  
  
    \def\GL{\mathrm{GL}}    \def\dd{\displaystyle}   \def\spanning{\textnormal{-span}}
  \def\cB{\mathcal{B}}
 \def\veps{\varepsilon} 
  
\def\diag{\mathrm{diag}}   \def\ss{\scriptstyle}
\newcommand{\cT}{\mathcal{T}}  

\newcommand{\corank}{\mathrm{corank}}
\newcommand{\rank}{\mathrm{rank}}

\newcommand{\cG}{\mathcal{G}}

\newcommand{\fkn}{\mathfrak{n}}
\newcommand{\supp}{\mathrm{supp}}

\newcommand{\Tr}{\mathrm{Tr}}
\newcommand{\cP}{\mathcal{P}}
\newcommand{\cM}{\mathcal{M}}

\newcommand{\Null}{\mathrm{Null}}
\newcommand{\ann}{\mathrm{ann}}

\makeatletter
\renewcommand{\@makefnmark}{\mbox{\textsuperscript{}}}
\makeatother

\allowdisplaybreaks[1]

\UseCrayolaColors

\begin{document}
\title{Supercharacter formulas for pattern groups}
\author{Persi Diaconis and Nathaniel Thiem}
\date{}

\maketitle

\begin{abstract}
C. Andre and N. Yan introduced the idea of a supercharacter theory to give a tractable substitute for character theory in wild groups such as the unipotent uppertriangular group $U_n(\FF_q)$.  In this theory superclasses are certain unions of conjugacy classes, and supercharacters are a set of characters which are constant on superclasses.   This paper gives a character formula for a supercharacter evaluated at a superclass for pattern groups and more generally for algebra groups.
\end{abstract}

\section{Introduction\protect\footnote{MSC 2000: 20C99 (05E)}\protect\footnote{Keywords: supercharacters, superclasses, finite unipotent group, algebra group,posets}}

Let $U_n(\FF_q)$ be the group of uppertriangular $n\times n$ matrices with entries in the finite field $\FF_q$, and with ones on the diagonal.  While describing the conjugacy classes or irreducible characters of $U_n(\FF_q)$ is a well-known wild problem, Carlos Andr\'e \cite{An95,An99,An02} and Ning Yan \cite{Ya01,Ya06} have shown that certain unions of conjugacy classes (here-after superclasses) and certain characters (here-after supercharacters) have an elegant theory that is rich enough to handle some Fourier analysis problems classically needing the full character table and yet tractable enough to admit a closed form formula for a supercharacter at a superclass (see also \cite{ADS04}).

Diaconis and Isaacs \cite{DI06} abstracted supercharacter theory to algebra groups, a class of groups of the form $\{1+X\ \mid\ X\in \fkn\}$, where $\fkn$ is a nilpotent $\FF_q$-algebra.   In the resulting theory, 
restriction of supercharacters and tensor products of supercharacters decompose as nonnegative integer combinations of supercharacters.  Furthermore, there is a notion of superinduction that is dual to restriction of supercharacters.  However, instead of giving a nice formula as in the case of $U_n(\FF_q)$, they leave the character formula in the form of an orbit sum.

The present paper combines the above two perspectives by viewing algebra groups as subgroups of $U_n(\FF_q)$.  The primary focus of this paper is on a large class of subgroups called pattern subgroups.  Let $J\subseteq \{(i,j)\ \mid\ 1\leq i<j\leq n\}$  be a subset that is \emph{closed} in the sense that $(i,j),(j,k)\in J$ implies $(i,k)\in J$ (ie. $J$ is a partial order on $\{1,2,\ldots, n\}$).  The \emph{pattern group} $U_J$ is the subgroup of $U_n(\FF_q)$ consisting of matrices whose $(i,j)$th entry is permitted to be nonzero only if $(i,j)\in J$.  In fact, pattern groups are type $A$ root subgroups, and C. Andr\'e  and A. Neto have constructed a supercharacter theory for types $B$, $C$, and $D$ \cite{AN06}.

We give a reasonably explicit formula for a supercharacter evaluated at a superclass for an arbitrary algebra group.  It is, however, easiest to state for the special case of pattern groups.  Both supercharacters and superclasses can be indexed by certain (different) subsets of $U_J$.   For $x,y,\in U_J$, let $\chi^y$ be the supercharacter corresponding to $y$; then there exists an explicit $|J|\times|J|$ matrix $M$ and vectors $a,b\in \FF_q^{|J|}$ depending on $x$ and $y$ such that
\begin{equation}
\overline{\chi^y(x)}=\left\{\begin{array}{ll} \dd\frac{\chi^y(1)}{q^{\rank(M)}}\prod_{(i,j)} \theta(b_{ij}'b_{ij}+x_{ij}y_{ij}), & \text{if $Mb'=a$ and the product is well-defined,}\\ 0, & \text{otherwise,}\end{array}\right.
\label{CharacterFormulaForm}
\end{equation}
where $\theta:\FF_q\rightarrow \CC^\times$ is an isomorphism.    See Theorem 
\ref{PatternGroupCharacterFormula} and Theorem \ref{AlgebraGroupCharacterFormula} for explicit statements.  This formula reduces to the Andre/Yan formula in the case when $J=\{(i,j)\ \mid\ 1\leq i<j\leq n\}$, and to the usual character formula for the Heisenberg group in the case when $J=\{(i,j)\ \mid\ i=1\text{ or } j=n\}$.

Section 2 gives background on pattern groups, algebra groups and supercharacter theory.  We also characterize which subgroups of $U_n(\FF_q)$ are algebra groups.  Section 3 sets out our notation for pattern groups using the language of closed sets of roots.  Section 4 describes the orbits and coorbits of an algebra group on its radical $\fkn$, and dual $\fkn^*$.  This amounts to a careful study of row and column operations.

The main results are in Section 5, giving a formula in the form (\ref{CharacterFormulaForm}) for pattern groups.  If the partial order given by $J$ has no chains of length four, then we show that superclasses and supercharacters are conjugacy classes and irreducible characters.  In particular, our formulas give the classical characters of the extra special $p$-groups and certain unipotent radicals of Levi subgroups.  At the other extreme, our formulas reduce to the Andre-Carter formulae for $U_n(\FF_q)$.  Several other examples are also given explicitly.

Section 6 shows that the character formula also holds for general algebra subgroups of pattern groups (and thus for all algebra groups).

\subsection*{Acknowledgements}  We thank  Bob Guralnick, Marty Isaacs, and Ning Yan for their help, as well as Eric Marberg and Vidya Venkateswaran for their contributions to examples.

\section{Background}

This section gives an overview and pointers to the literature on three fundamental topics: pattern groups, algebra groups and supercharacters.

\subsection{Pattern groups}

Let 
$$U_n(\FF_q)=\{n\times n\text{ unipotent uppertriangular matrices with entries in $\FF_q$}\}.$$

A \emph{closed} subset of pairs
$$J\subseteq \{ (i,j)\ \mid\ 1\leq i<j\leq n\}$$
is a subset with the property that $(i,j),(j,k)\in J$ implies $(i,k)\in J$.  For $t\in \FF_q$, let $x_{ij}(t)\in U_n(\FF_q)$ denote the matrix with ones on the diagonal, $t$ in the $(i,j)$th position and zeroes everywhere else.  Then for any closed subset $J\subseteq \{ (i,j)\ \mid\ 1\leq i<j\leq n\}$, the pattern group $U_J$ is
\begin{align*}
U_J&=\{(u_{ij})\in U_n(\FF_q)\ \mid\ (i,j)\notin J\text{ implies }u_{ij}=0\}\\
&=\langle x_{ij}(t)\ \mid\ t\in \FF_q, (i,j)\in J\rangle.
\end{align*} 
The closedness of $J$ implies that $U_J$ is a subgroup of $U_n(\FF_q)$.  It is also clear that
$$|U_J|=q^{|J|}.$$

Note that the set $J$ gives rise to a poset $\cP$ on $\{1,2,\ldots, n\}$, where $i\leq_{\cP} j$ if and only if $(i,j)\in J$.   The transitivity property of $\cP$ is the same as the closed property of $J$.  In the language of posets, these pattern groups are similar to Rota's incidence algebras \cite{SO97}.  

Alternatively, if $R^+=\{\veps_i-\veps_j\ \mid\ 1\leq i< j\leq n\}$ are the usual positive roots of type $A$.  Then $J\subseteq R^+$ is closed if and only if $\alpha,\beta\in J$ and $\alpha+\beta\in R^+$ implies $\alpha+\beta\in J$.  From this point of view $U_J$ is the subgroup of $U_n(\FF_q)$ generated by the one-parameter subgroups $\langle x_\alpha(t) \ \mid\ t\in \FF_q\rangle$ corresponding to $\alpha\in J$.  

\vspace{.25cm}

\noindent\textbf{Examples.} Many naturally occurring subgroups are pattern groups.
\begin{enumerate}
\item If $J=\{(i,j)\ \mid\ 1\leq i<j\leq n\}$, then the corresponding poset is
$$\cP=
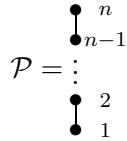
\xy<0cm,.8cm>\xymatrix@R=.4cm@C=.4cm{*={\bullet}\ar @{-} [d] & *={\ss n}\\ *={\bullet} \ar @{} [dd]|(.4)*{\vdots} & *={\ss n-1}\\ *={} & *={}\\    *={\bullet} \ar@{-} [d] & *={\ss 2}\\ *={\bullet} & *={\ss 1}}\endxy$$
and $U_J=U_n(\FF_q)$.
\item  The commutator subgroup $U_n'\subseteq U_n(\FF_q)$ is equal to the Frattini subgroup $\Phi(U_n)\subseteq U_n(\FF_q)$.  In this case, $U_n'=\Phi(U_n)=U_J$ where 
$$J=\{(i,j)\ \mid\ 1\leq i<i+1< j\leq n\},\qquad \text{or}\qquad \cP=\ 
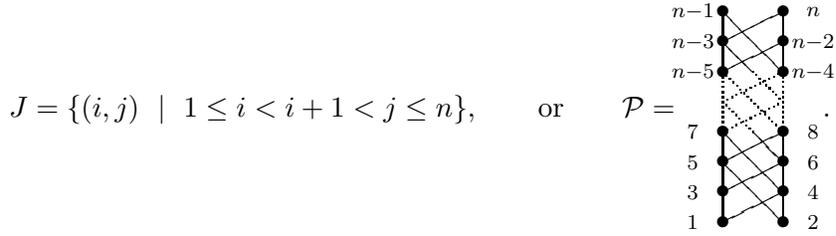
\xy<0cm,1.3cm> \xymatrix@R=.4cm@C=.4cm{*={\ss n-1} & *={\bullet} & *={} & *={\bullet} & *={\ss n}\\ 
*={\ss n-3} & *={\bullet} \ar @{-} [urr]  & *={} & *={\bullet} & *={\ss n-2}\\
*={\ss n-5} & *={\bullet}  \ar @{-} [urr]  \ar @{-} [uu]  & *={} \ar @{-} [ul] & *={\bullet}  \ar @{-} [uull]  \ar @{-} [uu]  & *={\ss n-4}\\
*={} & *={} \ar @{.} [urr] & *={} & *={} \ar @{.} [ul] & *={}\\
*={\ss 7} & *={\bullet} \ar @{.} [urr] \ar @{.} [uu] & *={} \ar @{.} [ul] & *={\bullet} \ar@{.} [uull] \ar @{.} [uu] & *={\ss 8}\\
*={\ss 5} & *={\bullet} \ar @{-} [urr] & *={} & *={\bullet} \ar @{-} [ul] & *={\ss 6}\\
*={\ss 3} & *={\bullet} \ar @{-} [urr]  & *={} & *={\bullet}\ar @{-} [uull]  & *={\ss 4} \\
*={\ss 1} & *={\bullet}\ar @{-} [uuu]\ar @{-} [urr] & *={} & *={\bullet} \ar @{-} [uull] \ar @{-} [uuu] & *={\ss 2}}\endxy\ .
$$
In general, if $J$ is closed, then $U_J'=\Phi(U_J)=U_{J'}$  \cite{DS93}, where
$$J'=\{(i,k)\ \mid\ (i,j),(j,k)\in J, \text{ for some $j$}\}.$$
Let $\{t_1,t_2,\ldots, t_r\}$ be a basis for $\FF_q$ as an $\FF_p$-vector space (where $q=p^r$).  From the description of the Frattini subgroups, a minimal (Frattini) generating set for $U_J$ may be chosen as
$$\{x_{ij}(t)\ \mid\ (i,j)\in J, t\in \{t_1,t_2,\ldots,t_r\}, i\leq k\leq j\text{ implies $(i,k)\notin J$ or $(k,j)\notin J$}\}.$$
Note that the set of $(i,j)$ that come up in these generators is the set of covering relations of the corresponding poset $\cP$.
\item  The center of $U_n(\FF_q)$ is $U_J$ for $J=\{(1,n)\}$, or 
$$\cP=
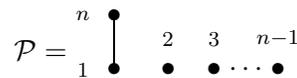
\xy<0cm,.5cm> \xymatrix@R=.2cm@C=.2cm{*={\ss n} & *={\bullet} \ar @{-} [dd] \\ &  & & *={\ss 2} & & *={\ss 3} & & *={\ss n-1}\\ *={\ss 1} & *={\bullet} & & *={\bullet} & & *={\bullet} & *{\cdots} & *={\bullet}}\endxy$$
\item The upper and lower central series of $U_n(\FF_q)$ has terms given by the pattern groups
$$J_k=\{(i,j)\ \mid\ j-i\geq k\}.$$
\item In $U_n(\FF_q)$, Vera-Lopez and Arregi \cite{AV04} used a sequence subgroups given by 
$$G_{ij}=\langle x_{kl}(t)\ \mid\ t\in \FF_q, i\leq k, j\leq l, (i,j)\neq (k,l)\rangle.$$
In this case,
$$1=G_{1n}\triangleleft G_{1,n-1}\triangleleft\cdots \triangleleft G_{12}\triangleleft G_{2n}\triangleleft\cdots\triangleleft G_{n-1,n}\triangleleft U_n(\FF_q)$$
gives a central series for $U_n(\FF_q)$ with each factor isomorphic to $\FF_q^+$.  Furthermore, every subgroup is a pattern group.
\end{enumerate}

The following elegant characterization of pattern groups was communicated to us by Bob Guralnick.

\begin{proposition}\label{PatternSubgroupsAsSubgroups}
For $q\geq 3$, a subgroup $U$ of $U_n(\FF_q)$ is a pattern group if and only if $U$ is invariant under conjugation by diagonal matrices $T$ in the general linear group $GL_n(\FF_q)$.
\end{proposition}
\begin{proof}
Note that
$$\diag(t_1,t_2,\ldots, t_n)(u_{ij})\diag(t_1^{-1},t_2^{-1},\ldots,t_n^{-1})=(t_iu_{ij}t_j^{-1}).$$
Thus, if $U$ is a pattern group, then it is invariant under the action of $T$.

Suppose $U$ is invariant under $T$.  Let $u=(u_{ij})\in U$ be nontrivial.  It suffices to show that for every $u_{ij}\neq 0$, the group $\langle x_{ij}(t)\ \mid\ t\in \FF_q\rangle$ is in $U$.  Let $i$ be minimal so that $u_{ij}\neq 0$ for some $i<j\leq n$.  Let
$$h_i(t)=\diag(\underbrace{1,\ldots, 1}_{i-1}, t, 1, \ldots, 1).$$ 
Then $h_i(t) u h_i(t^{-1})\in U$ has the effect of multiplying the $i$th row of $u$ by $t$.  In particular, if $t\neq 1$, then 
$$u'=h_i(t)uh_i(t^{-1})u^{-1}\in U\qquad \text{satisfies} \qquad \text{for $j<k$, } u_{jk}'=0\quad \text{unless $j=i$.}$$ 
Now let $j$ be minimal  such that $u_{ij}'\neq 0$ (such a $j$ must exist by our choice of $i$).  Then $h_j(t_2) u' h_j(t_2^{-1})\in U$ has the effect of multiplying the $(i,j)$th entry by $t_2^{-1}$.  Thus, if $t_2\neq 1$,
$$u''=h_j(t_2) u' h_j(t_2^{-1})u'^{-1}\in U\qquad \text{satisfies} \qquad u_{kl}''=0\quad \text{unless $k=i$ and $l=j$.}$$
Thus, $u''=x_{ij}(t)$ for some $t\in \FF_q^\times$.  Since $U$ is invariant under $T$, 
$$\langle h_i((t_3)x_{ij}(t)h_i(t_3^{-1})\ \mid\ t_3\in \FF_q^\times\rangle=\langle x_{ij}(t_3t)\ \mid\ t_3\in \FF_q^\times\rangle=\langle x_{ij}(t)\ \mid\ t\in \FF_q\rangle\subseteq U.$$
Note that $u'''=ux_{ij}(-u_{ij})$ satisfies $u_{ij}'''=0$ and $u_{kl}'''=u_{kl}$ for all $(k,l)\neq (i,j)$.  We may therefore proceed inductively to find all the one-parameter subgroups in $U$.
\end{proof}
 
\noindent\textbf{Remarks.}
\begin{enumerate}
\item  The proof fails for $q=2$ because the proof requires $|\FF_q^\times|>1$.  In fact, if $q=2$ the diagonal condition is empty, and the proposition is false.   For example,
$$\left\{\left(\begin{array}{ccc} 1 & s & t\\ 0 & 1 & s\\ 0 & 0 & 1\end{array}\right)\ \mid\ s,t\in \FF_2\right\}\subseteq U_3(\FF_2),$$
is a subgroup but not a pattern group since it does not contain $x_{12}(1)$. 
\item Proposition \ref{PatternSubgroupsAsSubgroups} implies all characteristic subgroups of $U_n(\FF_q)$ are pattern subgroups.
\item The proof does not depend on the field (other than the size constraint).
 \end{enumerate}

\subsection{Algebra groups}

Algebra groups generalize pattern groups.  Let $\fkn$ be a finite dimensional nilpotent $\FF_q$-algebra.     The corresponding \emph{algebra group} $U_{\fkn}$ is the group
$$U_{\fkn}=\{1+X\ \mid\ X\in \fkn\} \quad \text{with multiplication given by} \quad (1+X)(1+Y)=1+X+Y+XY.$$  
Algebra groups have $q$-power order, and in \cite{Rob98} Robinson studied the number of conjugacy classes of an algebra group.  It is obvious that the center of an algebra group is an algebra group since $Z(U_{\fkn})=1+Z(\fkn)$.  We do not know if the commutator or Frattini subgroups of an algebra group are algebra groups.  The classification of nilpotent algebras is an impossible task, but many classes of examples are known (see Pierce \cite{Pi82}).  In this section we show that algebra groups are poset groups.  

By Engels theorem we may view $\fkn$ as a subalgebra of the set of $n\times n$ uppertriangular matrices with zeroes on the diagonal (for some suitably chosen $n$), so every algebra group is isomorphic to a subgroup of $U_n(\FF_q)$.  The following proposition characterizes which subgroups of $U_n(\FF_q)$ are algebra groups.  Indeed, we determine which subgroups of $U_J$ are algebra groups for general closed sets $J$.

Let $J^*=\{\phi:J\rightarrow \FF_q\}$.  For $\phi\in J^*$, let
\begin{equation*}
x_\phi  =\left(\begin{array}{ccc} 1 & & \phi(i,j)\\ & \ddots &\\ 0 & & 1\end{array}\right) \in U_J\qquad\text{and}\qquad
X_\phi =\left(\begin{array}{ccc} 0 & & \phi(i,j)\\ & \ddots &\\ 0 & & 0\end{array}\right)\in \fkn_J=U_J-1.
\end{equation*}

\begin{proposition}\label{AlgebraGroupCharacterization}  For closed $J$, 
let $H\subseteq U_J$ be a subgroup.  Then
$$V=\{\phi\in J^*\ \mid\ x_\phi\in H\}$$
 is an $\FF_q$-vector space if and only if $H$ is a sub-algebra group of $U_J$.
\end{proposition}
\begin{proof}
Suppose $V$ is a vector space. It suffices to show
$$\fkn_H=\{X_\phi\ \mid\ \phi\in V\}$$
is an $\FF_q$-algebra.  Note that
$$aX_\phi+bX_\rho\in \fkn_H$$
is equivalent to $V$ being a vector space.  Consider
$$X_\phi X_\rho=(x_\phi-1)(x_\rho-1)=x_\phi x_\rho-x_\phi-x_\rho+1.$$
Since $x_\phi x_\rho\in H_J$, $x_\phi x_\rho=x_\eta$ for some $\eta\in V$.  Thus,
$$X_\phi X_\rho=X_{\eta-\phi-\rho}\in \fkn_H.$$

Suppose $H_J$ is a sub algebra group of $U_J$.  If $t_\rho, t_\mu\in V$ and $a,b\in \FF_q$, then 
$$1+(aX_\rho+bX_\mu)\in H_J$$
implies that $a\rho+b\mu\in V$.
\end{proof}

\noindent\textbf{Remarks.}
\begin{enumerate}
\item  Suppose $H\subseteq U_n(\FF_q)$ is an algebra group.  Then there is a natural poset $\cP$ on $\{1,2,\ldots, n\}$ given by $i\leq_\cP k$ if either
\begin{enumerate}
\item[(a)] there exists $x_\phi\in H$ such that $\phi_{ik}\neq 0$,
\item[(b)] there exist $x_\phi,x_\rho\in H$ such that $\phi_{ij}\neq 0$ and $\rho_{jk}\neq 0$ for some $j\in \{i+1,i+2,\ldots, k-1\}$.
\end{enumerate}
If $J=\{(i,j)\in \cP\times\cP\ \mid\ i<_\cP j\}$, then $U_J$ is the smallest pattern subgroup of $U_n(\FF_q)$ that contains $H$.  For example, if $H$ is the algebra subgroup
$$H=\left\{\left(\begin{array}{cccc} 1 & a & a & 0\\ 0 & 1 & 0 & -a\\ 0 & 0 & 1 & a\\ 0 & 0 & 0 & 1\end{array}\right)\ \mid a\in \FF_q\right\}\subseteq U_4(\FF_q),$$
then the corresponding poset is
$$\cP=\xy<0cm,.75cm>\xymatrix@R=.25cm@C=.25cm{*={} & {4}\\ {2} \ar @{-} [ur] & *={} & {3}\ar @{-} [ul]\\ *={} & {1} \ar @{-} [ul] \ar @{-} [ur]}\endxy\quad\text{with pattern group}\quad U_J=\left\{\left(\begin{array}{cccc} 1 & a & c & e\\ 0 & 1 & 0 & d\\ 0 & 0 & 1 & b\\ 0 & 0 & 0 & 1\end{array}\right)\ \mid a,b,c,d,e\in \FF_q\right\}.$$
\item Not every $p$-group is an algebra group.  In fact, one can use Proposition \ref{AlgebraGroupCharacterization} to show that the maximal unipotent subgroup of ${\rm Sp}_{2n}(\FF_q)$ is not an algebra group (for $n>1$).  In this way, algebra groups are a ``type A" phenomenon.
\item Much of the theory can be adapted when we replace $\FF_q$ by a finite radical ring, as explored in \cite{AnNi06}.
\end{enumerate}
\subsection{Supercharacters for algebra groups}\label{SectionIntroSuperCharacters}

Determining conjugacy classes and characters of $p$-groups contains intractable obstacles.  For $U_n(\FF_q)$, Andr\'e \cite{An95,An99,An02}, Yan \cite{Ya01} and Arias-Castro et al \cite{ADS04} found that certain unions of conjugacy classes and corresponding sums of irreducible characters give a tractable, useful theory.  This was abstracted to algebra groups in \cite{DI06}.  We give a brief synopsis.

First suppose that $G$ is an arbitrary finite group.   A \emph{supercharacter theory} for $G$ is a partition $\kappa$ of the conjugacy classes, and a partition $\kappa^\vee$ of the irreducible characters such that
\begin{enumerate}
\item[(a)] The identity element is in its own block in $\kappa$,
\item[(b)] $|\kappa|=|\kappa^\vee|$,
\item[(c)] For each block $K\in \kappa^\vee$, there exists a corresponding character $\chi^K$ which is a positive linear combination of the characters in $K$ such that $\chi^K$ is constant on  the blocks of $\kappa$.
\end{enumerate}
There are many general examples of supercharacter theories, including the partitioning of $G$ and its characters under the action of a group of automorphisms of $G$.   The following specific construction for algebra groups specializes to the construction of Andr\'e-Yan(see \cite{DI06}).  

If $U_{\fkn}$ is an algebra group with corresponding nilpotent algebra $\fkn$, then $U_\fkn\times U_\fkn$ acts on $\fkn$ by left and right multiplication.  A \emph{superclass} of $U_\fkn$ is a subset $1+O$, where $O$ is a two-sided orbit in $\fkn$.  Note that every superclass is in fact a union of conjugacy classes and $1$ is in its own superclass.  

To define supercharacters, note that $U_\fkn\times U_\fkn$  acts on the space of linear functionals $\fkn^*$ by
$$(u \lambda v)(X)=\lambda(u^{-1}Xv^{-1}), \qquad\text{where $u,v\in U_\fkn$, $\lambda\in \fkn^*$, and $X\in \fkn$}.$$
Fix a nontrivial homomorphism $\theta:\FF_q^+\rightarrow \CC^\times$.  For $\lambda\in \fkn^*$, define the supercharacter $\chi^\lambda$ to be
$$\chi^\lambda=\frac{|\lambda U_\fkn|}{|U_\fkn \lambda U_\fkn|}\sum_{\mu\in U_{\fkn}\lambda U_{\fkn}} \theta\circ \mu.$$
Then every irreducible character appears as a constituent of exactly one supercharacter (not obvious), and the supercharacters are constant on superclasses (follows from the definition).

The classical orbit method of Kirillov identifies conjugacy classes and irreducible characters with orbits of $U_\fkn$ acting on $\fkn$ by conjugation.  The problem with this construction is that describing the orbits is a provably wild problem.  Andr\'e and Yan have shown that the $U_\fkn\times U_\fkn$-orbits for $U_n(\FF_q)$ are indexed by labeled set partitions and that natural quantities such as dimensions and intertwining numbers are described in terms of elegant combinatorics (rivaling the tableaux combinatorics of the symmetric groups).   We hope that the following developments show that some of this carries over to general algebra groups.

Diaconis and Isaacs describe in \cite{DI06} that the above superclasses and supercharacters for algebra groups form a ``nice" supercharacter theory in the following sense
\begin{enumerate}
\item[(a)] The character of the regular representation is the sum of all the distinct supercharacters.  In particular, supercharacters are orthogonal with respect to the usual inner product,
\item[(b)] The restriction of a supercharacter to any algebra subgroup is a sum of supercharacters with nonnegative integer coefficients,
\item[(c)] There is a notion of ``superinduction" for which Frobenius reciprocity holds,
\item[(d)] The tensor product of supercharacters is a linear combination of supercharacters with nonnegative integer coefficients.
\end{enumerate}
In addition, \cite{DI06} give useful criteria for determining when a superclass or supercharacter is in fact a conjugacy class or irreducible character.   We apply this in Corollary \ref{FullUCharacters} to give a simple necessary and sufficient condition for a supercharacter of $U_n(\FF_q)$ to be irreducible.

The theory is rich enough to permit analysis of natural problems.  For example, consider a pattern group $U_J$.  Fix a basis $\{t_1,\ldots, t_r\}$ of $\FF_q$ as an $\FF_p$-vector space ($q=p^r$).  Then the set 
$$\{x_{ij}(t)\ \mid\ t\in \{t_1,\ldots, t_r\}, (i,j)\in J, \text{there is no pair $(i,k),(k,j)\in J$}\}$$  
is a generating set for $U_J$.  It follows from \cite[Corollary 3.5] {DI06} that the conjugacy class containing  one of these generators is a superclass.  Thus, supercharacter theory can be used to explicitly diagonalize the conjugacy walk on $U_J$.  Then one may attempt comparison theory \cite{DS93} to analyze the original walk.   This program was carried out in \cite{ADS04} for $U_n(\FF_q)$.

\section{Pattern groups}

This section gives two useful orderings on pairs $(i,j)$ and defines pattern groups.

\subsection{A Lie theoretic perspective}

Let
$$R^+=\{\veps_i-\veps_j\ \mid\ 1\leq i<j\leq n\},$$
be a set of roots, which has a total order given by
\begin{equation} \label{LinearOrdering}
\veps_r-\veps_s\leq \veps_i-\veps_j, \qquad \text{if $r>i$ or if $r=i$, and $s>j$}.
\end{equation}  
If the roots are indexed by pairs $(i,j)$, then this ordering can be pictured as
$$(n-1,n)<(n-2,n)<(n-2,n-1)<\cdots <(1,n)<(1,n-1)<\cdots < (1,2).$$
\begin{lemma} \label{OrderingLemma} The total order $\leq$ satisfies
\begin{enumerate}
\item[(a)] If $\alpha<\beta$ and $\alpha+\beta\in R^+$, then $\alpha< \alpha+\beta< \beta$.
\item[(b)] If $\alpha<\beta$, $\alpha+\beta\in R^+$ and $\alpha+\beta+\gamma\in R^+$, then either $\gamma< \alpha$ or $\beta< \gamma$.
\end{enumerate}
\end{lemma}
\begin{proof} (a) The assumptions imply $\alpha=\veps_j-\veps_k$, $\beta=\veps_i-\veps_j$, so $\alpha+\beta=\veps_i-\veps_k$.  Since $j>i$, $\alpha<\alpha+\beta$ and since $k>j$, $\alpha+\beta< \beta$.  (b) The assumptions imply that $\alpha=\veps_j-\veps_k$, $\beta=\veps_i-\veps_j$, and either $\gamma=\veps_h-\veps_i$ or $\gamma=\veps_k-\veps_l$.  
\end{proof}

There is also a partial ordering on the roots given by
$$\alpha\prec \beta, \qquad \text{if $\beta-\alpha\in R^+$},$$
called the \emph{dominance order}.  This ordering can be pictured by hanging the uppertriangle from its righthand corner $(1,n)$.  entries just above the main diagonal are lowest, those on the second diagonal above with $(1,n)$ at the top of the Hasse diagram (see the figure below).  This partial ordering is not compatible with the linear ordering  \ref{LinearOrdering}, and should also not be confused with the underlying poset $\cP$ on the integers $\{1,2,\ldots, n\}$.

Let $\cG(R^+)$ be the  Haase diagram of the poset $\preceq$, which we can organize from smallest to greatest as we pass from left to right.  For example, if $n=4$, 
$$\cG(R^+)=
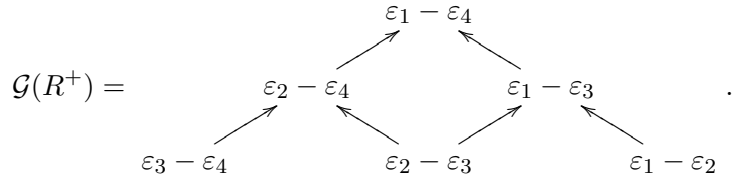
\xy<0cm,1cm>\xymatrix@R=.5cm@C=.25cm{ &  &  {\veps_1-\veps_4} & & \\ & {\veps_2-\veps_4}\ar  [ur] &  & {\veps_1-\veps_3} \ar   [ul]  & \\
{\veps_3-\veps_4} \ar  [ur] &  & {\veps_2-\veps_3} \ar  [ul] \ar  [ur]  & & {\veps_1-\veps_2}\ar  [ul] }\endxy.$$

A subset $J\subseteq R^+$ is \emph{closed} if $\alpha,\beta\in J$ implies that $\alpha+\beta\in J$.  Let $\cG(J)$ denote the subgraph of $\cG(R^+)$ with vertices $J$ and an edge from $\alpha$ to $\beta$ if $\beta-\alpha\in J$.  For example, when $n=4$, the subset $J=\{(1,2),(1,3),(1,4),(2,4),(3,4)\}$ is closed with  corresponding subgraph
$$\cG(J)=
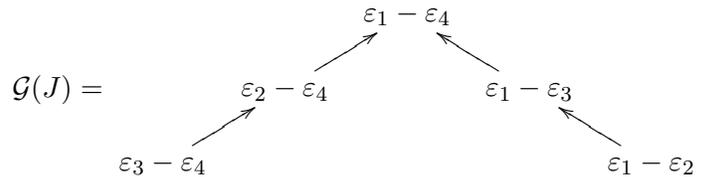
\xy<0cm,1cm>\xymatrix@R=.5cm@C=.25cm{ &  &  {\veps_1-\veps_4} & & \\ & {\veps_2-\veps_4}\ar  [ur] &  & {\veps_1-\veps_3} \ar   [ul]  & \\
{\veps_3-\veps_4} \ar  [ur] &  & *={} & & {\veps_1-\veps_2}\ar  [ul] }\endxy.$$

\subsection{Pattern groups}

For $\veps_i-\veps_j\in R^+$, let
$$X_{ij}=X_{\veps_i-\veps_j}=\text{ $n\times n$ matrix with 1 in the $(i,j)$th position and zeroes elsewhere}.$$
The nilpotent $\FF_q$-algebra 
$$\fkn_J=\FF_q\spanning\{X_\alpha\ \mid\ \alpha\in J\},$$
has relations
\begin{align}
X_\alpha^2&=0, & & \text{for $\alpha\in J$,} \label{LieSquare}\\
X_\alpha X_\beta & = 0, & & \text{for $\alpha,\beta\in J$ such that $\alpha< \beta$ OR $\alpha+\beta\notin R^+$,}\label{CorrectProduct}\\
X_\beta X_\alpha &  =  X_{\alpha+\beta} & & \text{for $\alpha,\beta,\alpha+\beta\in J$, $\alpha< \beta$ AND $\alpha+\beta\in R^+$.}\label{IncorrectProduct}
\end{align}

Let $J\subseteq R^+$ be closed, let $\FF_q$ be the finite field with $q$ elements, and let
$$J^*=\left\{\begin{array}{rcl}\phi:J &\rightarrow& \FF_q\\ \alpha & \mapsto & \phi_\alpha\end{array}\right\}.$$
Then
$$\fkn_J=\{X_\phi\ \mid\ \phi\in J^*\}\qquad\text{where}\qquad X_\phi=\sum_{\alpha\in J} \phi_\alpha X_\alpha.$$

The dual of $\fkn_J$,
$$\fkn_J^*=\{\lambda:\fkn_J\rightarrow \FF_q\ \mid\ \lambda\text{ $\FF_q$-linear}\},$$
has a basis $\{\lambda_\alpha: \fkn_J \rightarrow  \FF_q\ \mid\ \alpha\in J\}$ given by
$$\lambda_\alpha(X_\phi)=\phi_\alpha.$$
For $\eta\in J^*$, let
$$\lambda_\eta=\sum_{\alpha\in J} \eta_\alpha \lambda_\alpha,$$
so that $\fkn_J^*=\{\lambda_\eta\ \mid\ \eta\in J^*\}$.

For $\alpha\in J$ and $t\in \FF_q$, define
$$x_\alpha(t)=e^{tX_\alpha}=1+tX_\alpha.$$
Note that if we order the product by the total order $<$ (ie. begin with the smallest, multiply on the right by the next smallest, etc.), then relation (\ref{CorrectProduct}) implies
$$\prod_{\alpha\in J} x_\alpha(t_\alpha)=1+\sum_{\alpha\in J} t_\alpha X_\alpha,$$
so we may define the unipotent subgroup
\begin{align*}
U_J & = \langle x_\alpha(t) \ \mid\ \alpha\in J, t\in \FF_q\rangle\\
&= 1+\fkn_J.
\end{align*}
For $\phi\in J^*$, let
$$x_\phi=1+X_\phi=\prod_{\alpha\in J} x_\alpha(\phi_\alpha)\in U_J,$$
where again the product is ordered according to $<$.  Then
$$U_J=\{x_\phi\ \mid\ \phi\in J^*\}\qquad \text{and}\qquad \fkn_J=\{X_\phi\ \mid\ \phi\in J^*\}.$$
The generators of $U_J$ satisfy the commutation relation
\begin{equation}
[x_\alpha(a),x_\beta(b)]=\left\{\begin{array}{ll} x_{\alpha+\beta}(-ab), & \text{if $\alpha+\beta\in J$, $\alpha  < \beta$,}\\ x_{\alpha+\beta}(ab),  & \text{if $\alpha+\beta\in J$, $\alpha>\beta$}\\ 1, & \text{otherwise,}\end{array}\right.\label{ChevalleyCommutation}
\end{equation}
and the additive relation
\begin{equation} x_\alpha(a)x_\alpha(b)=x_\alpha(a+b). \label{RootSum}\end{equation}

\section{Pattern group orbits and co-orbits}

Both superclasses and supercharacters are indexed by two-sided orbits.  This section develops a description of these orbits.

\subsection{Pattern group orbits}\label{Orbits}

The group $U_J$ acts on $\fkn_J$ by left and right multiplication.  Let
$$O_\phi=\{\text{$U_J\times U_J$-orbit containing $X_\phi$}\}.$$

Of course, multiplying on the right by $x_\beta(t)$ where $\beta=(i,j)$ and $t\in \FF_q$ adds $t$ times column $i$ to column $j$ , while adding on  the left adds $t$ times row $j$ to row $i$.  In the present notation we obtain the following lemma.

\begin{lemma}[Row and column reducing]  \label{RowReducing} Let $\beta\in J$, $t\in \FF_q$, and $\phi\in J^*$.  Then
\begin{enumerate}
\item[(a)]  $X_\phi x_\beta(t)=X_{\phi'}$, where
$$\phi_\gamma'=\left\{\begin{array}{ll}\phi_\gamma+t\phi_{\gamma-\beta},  & \text{if $\gamma-\beta\in J$ AND $\beta< \gamma-\beta$,}\\ \phi_\gamma & \text{otherwise,}\end{array}\right.\qquad \gamma \in J,$$
\item[(b)] $x_\beta(t)X_\phi=X_{\phi'}$, where
$$\phi_\gamma' =\left\{\begin{array}{ll} \phi_\gamma+t\phi_{\gamma-\beta},  & \text{if $\gamma-\beta\in J$ AND $\gamma-\beta< \beta$,}\\
\phi_\gamma & \text{otherwise,} \end{array}\right.\qquad \gamma \in J.$$
\end{enumerate}
\end{lemma}
\begin{proof}  (a) Let $J=A\cup\{\beta\}\cup B$ be such that $A< \beta< B$, and let $\phi_A$ and $\phi_B$ be the restrictions of $\phi$ to $A$ and $B$, respectively.  Then
$$X_\phi x_\beta(t)=(x_\phi-1)x_\beta(t)=x_A(\phi_A)x_\beta(\phi(\beta))x_B(\phi_B)x_\beta(t)-x_\beta(t).$$
By relation (\ref{RootSum}), 
\begin{align*} X_\phi x_\beta(t)&= x_A(\phi_A)x_\beta(\phi(\beta)+t)[x_B(\phi_B),x_\beta(t)]-x_\beta(t)\\
&=x_A(\phi_A)x_\beta(\phi(\beta)+t)[x_B(\phi_B),x_\beta(t)]-1-tX_\beta\\
&= x_A(\phi_A)x_\beta(\phi(\beta))[x_B(\phi_B),x_\beta(t)]-1.
\end{align*}
By Lemma \ref{OrderingLemma} (a), if $\alpha\in B$ and $\alpha+\beta\in J$, then $\alpha+\beta\in B$.  By relation (\ref{ChevalleyCommutation}), as $x_\beta(t)$ moves left through $x_B(\phi_B)$ the only new terms that crop up are of the form $x_{\alpha+\beta}(t\phi_B(\alpha))$, where $\alpha\in B$ and $\alpha+\beta\in J$.  By Lemma  \ref{OrderingLemma} (b), as $x_{\alpha+\beta}(t\phi_B(\alpha))$ moves from right to left, it commutes with all other terms until it hits $x_{\alpha+\beta}(\phi_B(\alpha+\beta))$.   Thus, every $x_\gamma(\phi_B(\gamma))$ in $B$ such that $\gamma-\beta\in J$, becomes 
$x_\gamma(\phi_B(\gamma)+t\phi_B(\gamma-\beta))$, and all the other terms remain the same.  

The proof of (b) is similar.
\end{proof}

For the poset version of Lemma \ref{RowReducing}, write 
$$(i_1,i_2,\ldots, i_r)\in \cP\qquad\text{if}\qquad (i_1,i_2),(i_2,i_3),\ldots, (i_{r-1},i_r)\in J,$$
so that  $(i_1,i_2,\ldots, i_r)\in \cP$ if $(i_1,i_2,\ldots, i_r)$ is an $r$-chain in the poset.

\begin{lemma}
Let $x_{\veps_j-\veps_k}(t)\in U_J$ and $X_\phi\in \fkn_J$.  Then 
\begin{enumerate}
\item[(a)]  $X_\phi x_{\veps_j-\veps_k}(t)=X_{\phi'}$, where
$$\phi_{il}' =\left\{\begin{array}{ll} \phi_{il}+t\phi_{ij},  & \text{if $l=k$ AND $(i,j,l)\in \cP$,}\\
\phi_{il} & \text{otherwise,}\\ \end{array}\right.\qquad (i,l)\in \cP,$$
\item[(b)] $x_{\veps_j-\veps_k}(t)X_\phi=X_{\phi'}$, where
$$\phi_{il}' =\left\{\begin{array}{ll} \phi_{il}+t\phi_{kl},  & \text{if $i=j$ AND $(i,k,l)\in \cP$,}\\ \phi_{il} & \text{otherwise,}\end{array}\right.\qquad (i,l) \in \cP.$$
\end{enumerate}
\end{lemma}

By iterating Lemma \ref{RowReducing}, we obtain

\begin{theorem} \label{LeftRightAction} Let $\phi,\rho\in J^*$.  Then
\begin{enumerate}
\item[(a)]  $x_\rho X_\phi=X_{\phi'}$, where
$$\phi_{il}'= \phi_{il}+\sum_{(i,j,l)\in\cP} \rho_{ij}\phi_{jl},$$  
\item[(b)]  $X_\phi x_\rho= X_{\phi'}$ where
$$\phi_{il}'= \phi_{il}+\sum_{(i,j,l)\in\cP} \phi_{ij}\rho_{jl}.$$
\end{enumerate}
\end{theorem}
\begin{proof}
(a)  List the roots in $J$ according to the total order so 
$$\beta_{|J|}\leq \cdots \leq \beta_2\leq \beta_1.$$
 We can apply Lemma \ref{RowReducing} iteratively
 $$\left(x_{\beta_{|J|}}(\rho(\beta_{|J|})) \biggl(\cdots\biggl( x_{\beta_2}(\rho(\beta_2))\biggl(x_{\beta_1}(\rho(\beta_1))X_\phi\biggr)\biggr)\cdots\biggr)\right). $$
Induct on $1\leq r\leq |J|$.  By induction, if
$$X_{\phi^{(r-1)}}=\left(x_{\beta_{r-1}}(\rho(\beta_{r-1})) \biggl(\cdots\biggl( x_{\beta_2}(\rho(\beta_2))\biggl(x_{\beta_1}(\rho(\beta_1))X_\phi\biggr)\biggr)\cdots\biggr)\right),$$
then
$$\phi^{(r-1)}_{il}=\phi_{il} + \sum_{(i,j,l)\in \cP\atop \beta_r<(i,j)} \rho_{ij} \phi_{jl}.$$
Let $X_{\phi^{(r)}}=x_{\beta_r}(\rho(\beta_r)) X_{\phi^{(r-1)}}$.  By Lemma \ref{RowReducing}, 
$$\phi^{(r)}_{il} = \left\{\begin{array}{ll} (\phi^{(r-1)}_{il} +\rho_{ab}\phi^{(r-1)}_{bl}, & \text{if $a=i$ and $(a,b,l)\in \cP$,}\\ \phi^{(r-1)}_{il}, & \text{otherwise}.\end{array}\right.$$
However, if $a=i$ and $(a,b,l)\in \cP$, then by the choice of ordering, the set $\{(b,k,l)\in \cP\ \mid\  \beta_r<(b,k)\}$ must be empty (since $a\leq b$).  Thus,
\begin{align*}
\phi^{(r)}_{il} &= \left\{\begin{array}{ll} \phi^{(r-1)}_{il} +\rho_{ab}\phi_{bl}, & \text{if $a=i$ and $(a,b,l)\in \cP$,}\\ \phi^{(r-1)}_{il}, & \text{otherwise},\end{array}\right.\\
&= \phi_{il}+\sum_{(i,j,l)\in \cP\atop \beta_r\leq (i,j)} \rho_{ij}\phi_{jl},
\end{align*}
as desired.

The proof for (b) is similar.
\end{proof}
 
 By first applying (a) and the (b) of the Theorem to $X_\phi$, we obtain the following corollary.
 
\begin{corollary} \label{CompleteLeftRightAction}
Let $\phi, \tau, \rho\in J^*$.  If $X_{\phi'}=x_\tau X_\phi x_{\rho}$, then 
$$\phi_{il}'=\phi_{il}+\sum_{(i,j,l)\in\cP} \tau_{ij}\phi_{jl}+\sum_{(i,k,l)\in\cP}\phi_{ik} \rho_{kl}+\sum_{(i,j,k,l)\in\cP} \tau_{ij}\phi_{jk}\rho_{kl} .$$
\end{corollary}

Define the matrices $M_\phi^L$ and $M_\phi^R$ by
\begin{align*}
(M_\phi^L)_{(i,l),(j,k)}  & = \left\{\begin{array}{ll} \phi_{kl}, & \text{if $i=j$, $(i,k,l)\in \cP$},\\ 0, & \text{otherwise,}\end{array}\right.  & & \text{for $(i,l),(j,k)\in J$,}\\
(M_\phi^R)_{(i,l),(j,k)}  & = \left\{\begin{array}{ll} \phi_{ij}, & \text{if $k=l$, $(i,j,l)\in \cP$},\\ 0, & \text{otherwise,}\end{array}\right.  & & \text{for $(i,l),(j,k)\in J$.}
\end{align*}
\begin{corollary}  Let $\phi\in J^*$.  Then
\begin{enumerate}
\item[(a)] The size of the right orbit containing $X_\phi$ is $q^{\rank(M_\phi^R)}$,
\item[(b)] The size of the left orbit containing $X_\phi$ is $q^{\rank(M_\phi^L)}$.
\end{enumerate}
\end{corollary}
\begin{proof}
(a)  Consider the vector space isomorphism
$$\begin{array}{rcl} v: \fkn_J & \longrightarrow & \FF_q^{|J|}\\ X_\rho & \mapsto & v_\rho,\end{array}\qquad \text{where}\qquad (v_\rho)_{(i,j)}=\rho_{ij}.$$
By Theorem \ref{LeftRightAction}, the equation
$$X_\phi x_\rho= X_\phi + X_\phi X_\rho$$
becomes
$$v(X_\phi x_\rho)=v_\phi+ M_\phi^R v_\rho.$$
Thus, the right orbit containing $X_\phi$ has the same size as
$$|\{M_\phi^R v\ \mid\ v\in \FF_q^{|J|}\}|=q^{\rank(M_\phi^R)}.$$
The proof of (b) is similar.
\end{proof}

\noindent\textbf{Cautionary example.}  For the full upper-triangular group, much  of the character theoretic information depends only on the ``shape" or the support of the particular $\phi$ or $\eta$ in $J^*$ \cite{ADS04}.  The following example shows that this does not hold for general pattern groups (at least in the obvious way).  Consider $q\geq 3$ and
$$
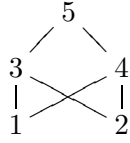
\xy<0cm,.7cm>\xymatrix@R=.3cm@C=.3cm{  & {5} &  \\ {3} \ar @{-} [ur] & & {4}\ar @{-} [ul] \\ {1} \ar @{-} [urr]  \ar @{-} [u]  & & {2} \ar @{-} [ull]\ar @{-} [u] }\endxy
$$  
The group $U_J$ is then the set of matrices of the form
$$\left(\begin{array}{ccccc} 1 & 0 & \ast & \ast & \ast \\ 0 & 1 & \ast & \ast  & \ast \\ 0 & 0 & 1 & 0  & \ast  \\ 0 & 0 & 0 & 1 & \ast \\ 0 & 0 & 0 & 0 & 1\end{array}\right).$$
The matrices
$$x_1 = \left(\begin{array}{ccccc} 1 & 0 & 1 & 1 & 0 \\ 0 & 1 & 1 & 1  & 0 \\ 0 & 0 & 1 & 0  & 0  \\ 0 & 0 & 0 & 1 & 0 \\ 0 & 0 & 0 & 0 & 1\end{array}\right) \qquad \text{and}\qquad x_2=\left(\begin{array}{ccccc} 1 & 0 & 2 & 1 & 0 \\ 0 & 1 & 1 & 1  & 0 \\ 0 & 0 & 1 & 0  & 0  \\ 0 & 0 & 0 & 1 & 0 \\ 0 & 0 & 0 & 0 & 1\end{array}\right)$$
have two-sided orbits given by
$$\left\{\left(\begin{array}{ccccc} 1 & 0 & 1 & 1 & a \\ 0 & 1 & 1 & 1  & a \\ 0 & 0 & 1 & 0  & 0  \\ 0 & 0 & 0 & 1 & 0 \\ 0 & 0 & 0 & 0 & 1\end{array}\right)\bigg|\ a\in \FF_q\right\}\quad\text{and}\quad\left\{\left(\begin{array}{ccccc} 1 & 0 & 2 & 1 & a \\ 0 & 1 & 1 & 1  & b \\ 0 & 0 & 1 & 0  & 0  \\ 0 & 0 & 0 & 1 & 0 \\ 0 & 0 & 0 & 0 & 1\end{array}\right)\bigg|\ a,b\in \FF_q\right\},$$
respectively, so although $x_1$ and $x_2$ have the same ``shape," they have different size superclasses.

\subsection{Pattern group co-orbits}

The group $U_J$ acts on the dual space $\fkn^*$ on the left and right by
$$x\lambda(X_\phi)y= \lambda(x^{-1}X_\phi y^{-1}), \qquad\text{for $x,y\in U_J$, $\lambda\in \fkn^*$.}$$
For $\eta\in J^*$, let
$$O^\eta=\{\text{$U_J\times U_J$-orbit containing $\lambda_\eta$}\}.$$

\begin{lemma}[Dual row and column reducing] \label{DualReducing} Let $\beta\in J$, $t\in \FF_q$, and $\eta:J\rightarrow \FF_q$.  Then
\begin{enumerate}
\item[(a)] $\lambda_\eta x_\beta(-t)= \lambda_{\eta^\prime}$, where
$$\eta_\alpha^\prime=\left\{\begin{array}{ll}\eta_\alpha+t\eta_{\alpha+\beta},  & \text{if $\alpha+\beta\in J$ AND $\beta< \alpha+\beta$,}\\ \eta_\alpha & \text{otherwise,} \end{array}\right.\qquad \alpha \in J,$$
\item[(b)] $x_\beta(-t)\lambda_\eta = \lambda_{\eta^\prime}$, where
$$\eta^\prime_\alpha=\left\{\begin{array}{ll}  \eta_\alpha+t\eta_{\alpha+\beta},  & \text{if $\alpha+\beta\in J$ AND $\alpha+\beta< \beta$,}\\ \eta_\alpha & \text{otherwise,}\end{array}\right.\qquad \alpha \in J.$$
\end{enumerate}
\end{lemma}
\begin{proof} (a) If $\alpha\leq \beta$ OR $\alpha+\beta\notin J$, then
$$\lambda_\eta(X_\alpha)x_\beta(-t)=\lambda_\eta(x_\alpha(1) x_\beta(t)-x_\beta(t))=\lambda_\eta(X_\alpha).$$
If $\beta< \alpha$ AND $\alpha+\beta\in J$, then
\begin{align*}
\lambda_\eta(X_\alpha)x_\beta(-t)&=\sum_{\gamma\in J}  \eta_\gamma\lambda_\gamma(X_\alpha x_\beta(t))\\
&= \sum_{\gamma\in J}  \eta_\gamma\lambda_\gamma(tX_{\alpha+\beta} + X_\alpha)\\
&= t\eta_{\alpha+\beta}\lambda_{\alpha+\beta}(X_{\alpha+\beta})+\eta_\alpha\lambda_\alpha(X_\alpha)\\
&=(\lambda_\eta+t\eta_{\alpha+\beta}\lambda_\alpha)(X_\alpha)
\end{align*}
Thus,
$$\lambda_\eta x_\beta(-t)=\lambda_{\eta^\prime}, \qquad \text{where} \qquad \eta^\prime_\alpha=\left\{\begin{array}{ll} \eta_\alpha+t\eta_{\alpha+\beta},  & \text{if $\alpha+\beta\in J$ AND $\beta< \alpha+\beta$,}\\ \eta_\alpha & \text{otherwise.}\end{array}\right.$$
The proof to (b) is similar.\end{proof}

The poset version of Lemma \ref{DualReducing} is
\begin{lemma} Let $x_{\veps_i-\veps_l}(t)\in U_J$, $t\in \FF_q$, and $\eta\in J^*$.  Then
\begin{enumerate}
\item[(a)] $x_{\veps_i-\veps_l}(-t)\lambda_\eta = \lambda_{\eta^\prime}$, where
$$\eta^\prime_{jk}=\left\{\begin{array}{ll}\eta_{jk}+t\eta_{ik},  & \text{if $l=j$ AND $(i,j,k)\in \cP$,}\\ \eta_{jk} & \text{otherwise,} \end{array}\right.\qquad (j,k)\in \cP,$$
\item[(b)] $\lambda_\eta x_{\veps_i-\veps_l}(-t)= \lambda_{\eta^\prime}$, where
$$\eta^\prime_{jk}=\left\{\begin{array}{ll}  \eta_{jk}+t\eta_{jl},  & \text{if $i=k$ AND $(j,k,l)\in \cP$,}\\ \eta_{jk} & \text{otherwise,}\end{array}\right.\qquad (j,k)\in \cP.$$
\end{enumerate}
\end{lemma}

We have a similar result to Theorem \ref{LeftRightAction}, but in this case we act by $x_\rho^{-1}$.
\begin{theorem} \label{DualLeftRightAction} Let $\eta,\rho\in J^*$.  Then
\begin{enumerate}
\item[(a)]  $x_\rho^{-1} \lambda_\eta=\lambda_{\eta'}$, where
$$\eta_{jk}'= \eta_{jk}+\sum_{(i,j,k)\in\cP} \rho_{ij}\eta_{ik},$$  
\item[(b)]  $\lambda_\eta x_\rho^{-1}= \lambda_{\eta'}$ where
$$\eta_{jk}'= \eta_{jk}+\sum_{(j,k,l)\in\cP} \rho_{kl}\eta_{jl}.$$
\end{enumerate}
\end{theorem}
\begin{proof}
The proof is the same as for Theorem \ref{LeftRightAction}, with the reversed ordering.
\end{proof}

Combine (a) and (b) from Theorem \ref{DualLeftRightAction}, to obtain

\begin{corollary}  \label{CompleteLeftRightCoAction}
Let $\eta,\tau,\rho\in J^*$.  If $\lambda_{\eta'}=x_\tau^{-1}\lambda_{\eta}x_{\rho}^{-1}$, then
$$\eta_{jk}'=\eta_{jk}+\sum_{(i,j,k)\in \cP} \tau_{ij}\eta_{ik}+\sum_{(j,k,l)\in \cP} \eta_{jl} \rho_{kl}+\sum_{(i,j,k,l)\in \cP} \tau_{ij} \eta_{il}\rho_{kl}.$$
\end{corollary} 

Define the $|J|\times |J|$ matrices $M_L^\eta$ and $M_R^\eta$ by
\begin{align*}
(M_L^\eta)_{(j,k),(i,l)}  & = \left\{\begin{array}{ll} \eta_{ik}, & \text{if $l=j$, $(i,j,k)\in \cP$},\\ 0, & \text{otherwise,}\end{array}\right.  & & \text{for $(j,k),(i,l))\in J$,}\\
(M_R^\eta)_{(j,k),(i,l)}  & = \left\{\begin{array}{ll} \eta_{jl}, & \text{if $k=i$, $(j,k,l)\in \cP$},\\ 0, & \text{otherwise,}\end{array}\right.  & & \text{for $(j,k),(i,l)\in J$.}
\end{align*}
Note that 
$$(M_L^\eta)_{(j,k),(i,j)}=\eta_{ik} \quad \text{if and only if} \quad (M_R^\eta)_{(i,j),(j,k)}=\eta_{ik}.$$
Thus, $\rank((M_L^\eta)=\rank(M_R^\eta)$.  Define the \emph{corank} of $\eta$ to be
$$\corank(\eta)=\rank(M_L^\eta)=\rank(M_R^\eta).$$

\begin{corollary}  Let $\eta\in J^*$.  Then
\begin{enumerate}
\item[(a)] The size of the right orbit containing $\lambda_\eta$ is $q^{\corank(\eta)}$,
\item[(b)] The size of the left orbit containing $\lambda_\eta$ is $q^{\corank(\eta)}$.
\end{enumerate}
\end{corollary}

\noindent\textbf{Cautionary example.}  Two co-orbits with the same support shape can also have different co-orbit sizes (cf. cautionary example in Section \ref{Orbits}).  A similar argument to the one in \ref{Orbits} gives that 
$$
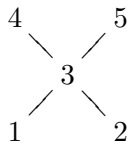
\xy<0cm,.7cm>\xymatrix@R=.3cm@C=.3cm{ {4} & & {5} \\ & {3}\ar @{-} [ul] \ar @{-} [ur] \\ {1} \ar @{-} [ur] & & {2} \ar @{-} [ul]}\endxy$$  
has co-orbits whose sizes depend on more than their shape.  Furthermore, while the number of orbits equals the number of co-orbits \cite[Lemma 4.1]{DI06}, the sizes of the two-sided orbits for an orbit or co-orbit indexed by the same symbol $\eta$ can differ (in fact, they usually will).

\section{A pattern group supercharacter formula}\label{SuperCharacterFormula}

This section states and proves our main theorem, a formula for a supercharacter on a superclass.  Following this, examples (Heisenberg and $U_n(\FF_q)$) show how the formula can be used.

The \emph{superclass} corresponding to $\phi\in J^*$ is
$$\{x_\rho\ \mid\ \rho\in O_\phi\},$$
and is a union of conjugacy classes in $U_J$. 

Fix a nontrivial homomorphism $\theta:\FF_q^+\rightarrow \CC^\times$.  For $\eta\in J^*$, let $\chi^\eta:U_J\rightarrow \CC$ be the map given by
\begin{align*}
\chi^\eta(x_\phi) & = \frac{q^{\corank(\eta)}}{|O^\eta|}\sum_{\mu\in O^\eta} \theta\bigl(\lambda_\mu(X_\phi)\bigr)\\
&= \frac{q^{\corank(\eta)}}{|O_\phi|}\sum_{\rho\in O_\phi} \overline{\theta\bigl(\lambda_\eta(X_\rho)}\bigr).
\end{align*}
By \cite[Theorems 5.6 and 5.8]{DI06}, these maps are the \emph{supercharacters} of $U_J$.  They are constant on superclasses, and are orthogonal under the usual inner product with the relation \cite[Lemma 5.9]{DI06},
$$\langle \chi^\eta,\chi^\mu\rangle=\delta_{\eta\mu} \frac{q^{2\corank(\eta)}}{|O^\eta|}.$$

To state the main theorem, we need one final piece of notation.  For $\phi,\eta\in J^*$, let $M_\phi^\eta$ be the $|J|\times |J|$ matrix given by
\begin{equation}\label{MeshMatrix}(M_\phi^\eta)_{(i,j)(k,l)}  = \left\{\begin{array}{ll} \phi_{jk}\eta_{il}, & \text{if $(i,j,k,l)\in \cP$,}\\ 0,  & \text{otherwise}.\end{array}\right.\end{equation}
Let $a_\phi^\eta\in \FF_q^{|J|}$ and $b_\phi^\eta\in \FF_q^{|J|}$  be given by
\begin{align}
(a_\phi^\eta)_{(i,j)} &= \sum_{(i,j,k)\in \cP} \phi_{jk}\eta_{ik} \label{aVector}\\
(b_\phi^\eta)_{(j,k)} &= \sum_{(i,j,k)\in \cP} \phi_{ij}\eta_{ik} \label{bVector}.
\end{align}
The functional $\phi$ \emph{meshes} with $\eta$ if
\begin{enumerate}
\item[(1)]  a solution to the equation $M_\phi^\eta x=-a_\phi^\eta$ exists,
\item[(2)]  $b_\phi^\eta$ is perpendicular to the nullspace of $M_\phi^\eta$ (with respect to the usual dot product).
\end{enumerate}
The following theorem gives a formula for a supercharacter on a superclass.

\begin{theorem}\label{PatternGroupCharacterFormula}
Let $\phi,\eta\in J^*$.  Let $M_\phi^\eta$, $a_\phi^\eta$, and $b_\phi^\eta$ be as in (\ref{MeshMatrix}), (\ref{aVector}), and (\ref{bVector}).  Then
$$\overline{\chi^\eta(x_\phi)}= \left\{\begin{array}{ll} \dd q^{\mathrm{corank}(\eta)-\mathrm{rank}(M_\phi^\eta)} \theta(b_0\cdot b_\phi^\eta) \hspace{-.5cm} \prod_{(i,j)\in \supp(\phi)\cap\supp(\eta)}  \hspace{-.5cm}\theta(\phi_{ij}\eta_{ij}),  & \text{if  $\phi$ meshes with $\eta$,}\\ 0, & \text{otherwise,}\end{array}\right.$$
where $b_0\in \FF_q^{|J|}$ satisfies $M_\phi^\eta b_0=-a_\phi^\eta$.
\end{theorem}
\begin{proof} 
Let $\FF_J=\FF_q^{|J|}$.
For $a,b\in \FF_J$, let $X_{ab}=X_{\phi'}$, where
$$\phi_{il}'=\phi_{il}+\sum_{(i,j,l)\in \cP} a_{ij} \phi_{jl} +\sum_{(i,k,l)\in \cP} \phi_{ik}b_{kl}+\sum_{(i,j,k,l)\in \cP}a_{ij} \phi_{jk}b_{kl}.$$
Since $|O_\phi|$ divides $q^{2|J|}$, we can overcount to get
$$\overline{\chi^\eta(x_\phi)}=\frac{q^{\mathrm{corank}(\eta)}}{|O_{\phi}|}\sum_{\rho\in \kappa_{\phi}}\theta(\lambda_\eta(X_\rho))
=\frac{q^{\mathrm{corank}(\eta)}}{q^{2|J|}}\sum_{a,b\in \FF_J}\theta(\lambda_\eta(X_{ab})).
$$
By the definition of $\lambda_\eta$ and $X_{ab}$
\begin{align*}
\overline{\chi^\eta(x_\phi)}
&=\frac{q^{\mathrm{corank}(\eta)}}{q^{2|J|}}\sum_{a,b\in \FF_J}\theta\left(\sum_{(i,l)\in J} \eta_{il}\biggl(\phi_{il} + \sum_{(i,j,l)\in \cP} a_{ij}\phi_{jl}+\phi_{ij}b_{jl}+\sum_{(i,j,k,l)\in \cP} a_{ij}\phi_{jk}b_{kl} \biggr)\right)\\
&=\frac{q^{\mathrm{corank}(\eta)}\theta_{\phi\eta}}{q^{2|J|}}\sum_{a,b\in \FF_J}\theta\left(\sum_{(i,l)\in J} \eta_{il}\biggl(\sum_{(i,j,l)\in \cP} a_{ij}\phi_{jl}+\phi_{ij}b_{jl}+\sum_{(i,j,k,l)\in \cP} a_{ij} \phi_{jk}b_{kl}\biggr)\right),
\end{align*}
where the second equality comes by collecting all summands that do not depend on $a$ and $b$,
$$\theta_{\phi\eta}=\theta\biggl(\sum_{(i,l)\in J} \eta_{il}\phi_{il}\biggr)=\prod_{(i,j)\in \supp(\phi)\cap\supp(\eta)} \theta(\eta_{il}\phi_{il}).$$
Collect summands that contain $a_{ij}$ for all $(i,j)\in J$,
\begin{align*}
\sum_{a,b\in \FF_J}&\theta\left(\sum_{(i,l)\in J} \eta_{il}\biggl(\sum_{(i,j,l)\in \cP} a_{ij}\phi_{jl}+\phi_{ij}b_{jl}+\sum_{(i,j,k,l)\in \cP} a_{ij} \phi_{jk}b_{kl}\biggr)\right)\\
&=\sum_{a,b\in \FF_J}\theta\left(\sum_{(i,j)\in J} a_{ij}\biggl(\sum_{(i,j,k)\in\cP} \phi_{jk}\eta_{ik} +\sum_{(i,j,k,l)\in \cP}  \phi_{jk}\eta_{il}b_{kl}\biggr)+\sum_{(j,k,l)\in\cP} \phi_{jk} \eta_{jl} b_{kl}\right)\\
&=\sum_{a,b\in \FF_J}\prod_{(i,j)\in J}\theta\left(a_{ij}\biggl(\sum_{(i,j,k)\in\cP} \phi_{jk}\eta_{ik} +\sum_{(i,j,k,l)\in \cP}  \phi_{jk}\eta_{il}b_{kl}\biggr)+ \sum_{(j,k,l)\in\cP} \phi_{jk} \eta_{jl} b_{kl}\right).
\end{align*}
Note that if for any $(i,j)\in J$, 
$$\sum_{(i,j,k)\in\cP} \phi_{jk}\eta_{ik} +\sum_{(i,j,k,l)\in \cP} \phi_{jk}\eta_{il}b_{kl} \neq 0,$$
then as we sum over all possible values of $a_{ij}\in \FF_q$, we sum over all roots of unity and
$$\sum_{a_{ij}\in \FF_q} \theta\left(a_{ij}\biggl(\sum_{(i,j,k)\in\cP} \phi_{jk}\eta_{ik} +\sum_{(i,j,k,l)\in \cP}  \phi_{jk}\eta_{il}b_{kl}\biggr)+ \sum_{(j,k,l)\in\cP} \phi_{jk} \eta_{jl} b_{kl}\right)=0.$$
If
$$\mathcal{S}=\{b\in\FF_J\ \mid\ M_\phi^\eta b=-a_\phi^\eta\},$$
then the equation
$$\sum_{(i,j,k)\in\cP} \phi_{jk}\eta_{ik} +\sum_{(i,j,k,l)\in \cP}  \phi_{jk}\eta_{il}b_{kl}=(a_\phi^\eta)_{ij}+(M_\phi^\eta b)_{ij},\qquad \text{for $(i,j)\in J$},$$
implies
\begin{equation*}
\overline{\chi^\eta(x_\phi)} = \frac{q^{\mathrm{corank}(\eta)}\theta_{\phi\eta}}{q^{2|J|}}\sum_{a\in \FF_J,b\in \mathcal{S}}\theta\left( \sum_{(j,k,l)\in\cP} \phi_{jk} \eta_{jl}b_{kl} \right)
= \frac{q^{\mathrm{corank}(\eta)}\theta_{\phi\eta}}{q^{|J|}}\sum_{b\in \mathcal{S}}\theta(b_\phi^\eta\cdot b).
\end{equation*}
Fix an element $b_0\in \mathcal{S}$.  Since every other vector in $\mathcal{S}$ is of the form $b_0+b'$, where $b'$ is in the nullspace of $M_\phi^\eta$, we have
$$\overline{\chi^\eta(x_\phi)} = \frac{q^{\mathrm{corank}(\eta)}\theta_{\phi\eta}\theta(b_0\cdot b_\phi^\eta)}{q^{|J|}}\sum_{b'\in \mathrm{Null}(M_\phi^\eta)}\theta(b'\cdot b_\phi^\eta).$$
Let $\{b_1',b_2',\ldots, b_r'\}$ be a basis for $\mathrm{Null}(M_\phi^\eta)$.  Then
\begin{equation*}
\sum_{b'\in \mathrm{Null}(M_\phi^\eta)}\theta(b'\cdot b_\phi^\eta)=\sum_{t\in \FF_q^r} \prod_{i=1}^r \theta(t_i b_i'\cdot b_\phi^\eta)= \prod_{i=1}^r \sum_{t_i\in \FF_q}\theta(t_i b_i'\cdot b_\phi^\eta).
\end{equation*}
Thus, if any $b_i'$ is not orthogonal to $b_\phi^\eta$, then the product is zero.  Thus, if $\phi$ meshes with $\eta$, we get
\begin{align*}
\overline{\chi^\eta(x_\phi)} &=  \frac{q^{\mathrm{corank}(\eta)}\theta_{\phi\eta}\theta(b_0\cdot b_\phi^\eta)}{q^{|J|}} |\mathrm{Null}(M_\phi^\eta)|= \frac{q^{\mathrm{corank}(\eta)}\theta_{\phi\eta} \theta(b_0\cdot b_\phi^\eta)}{q^{|J|}}q^{|J|-\mathrm{rank}(M_\phi^\eta)}\\
&=q^{\mathrm{corank}(\eta)-\mathrm{rank}(M_\phi^\eta)}\theta_{\phi\eta}\theta(b_0\cdot b_\phi^\eta),
\end{align*}
as desired.
\end{proof}

\noindent\textbf{Remarks.}
\begin{enumerate} 
\item In general, the matrices $M_\phi^\eta$, $M_\phi^R$, $M_R^\eta$, $M_\phi^L$, and $M_L^\eta$ have many zero rows and columns.  If one is careful, one can significantly reduce the dimensions of these matrices using the structure of the corresponding posets, but for expository purposes we omitted such ``simplifications."
\item  We have avoided choosing orbit and co-orbit representatives in stating Theorem \ref{PatternGroupCharacterFormula}.  In concrete cases, there appear to be natural choices (often with minimal support).  For example, the Heisenberg group, below, orbits and co-orbits can be identified with cosets of the center.  With $U_n(\FF_q)$, below, orbits and co-orbits can be identified with labeled set partitions as described in \cite{An95,An99,An01,An02,ADS04,Ya01}.
\end{enumerate}

\subsection{Examples}

\subsubsection*{Example 1: Heisenberg group.}  Let $H_n$ bet the group of order $q^{2n-3}$ represented as $n\times n$ upper triangular matrices with ones on the diagonal, entries in $\FF_q$, and non-zero entries allowed only in the top row or last column.   That is, let $H_n$ be the pattern group $U_J$ with
\begin{equation} \label{HeisenbergOrder} J=\{(n-1,n),(n-2,n),\ldots, (1,n), (1,n-1),\ldots, (1,2)\},\end{equation}
and corresponding poset
$$\cP= 
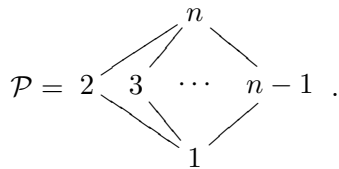
\xy<0cm,1cm>\xymatrix@R=.5cm@C=.25cm{*={} & *={} & {n}\\ {2} \ar @{-} [urr] & {3} \ar @{-} [ur] & {\cdots} & {n-1}\ar@{-} [ul]\\ *={} & *={} & {1} \ar @{-} [ull] \ar @{-} [ul] \ar @{-} [ur]} \endxy\ .$$
The structure and character theory of $H_n$  is well-known.  The center 
$$Z(H_n)= \text{the commutator of $H_n$} = \{ x_{1n}(t) \ \mid\  t\in \FF_q\}.$$
Thus, the conjugacy classes of $H_n$ consist of $Z(H_n)$ together with the $(q^{2n-4}-1)$ other cosets of $Z(H_n)$.    There are $q^{2n-4}$ linear characters and $(q-1)$ irreducible characters of degree $q^{n-2}$.  

Consider the supercharacter theory of $H_n$.  From Corollary \ref{CompleteLeftRightAction} (or direct computation), $x_\tau X_\phi x_{\rho}=X_{\phi'}$ with $\phi'_{il}=\phi_{il}$ unless $(i,l)=(1,n)$ and 
$$\phi_{1n}'=\phi_{1n}+\sum_{j=2}^{n-1} (\tau_{1j}\phi_{1n}+\phi_{1n}\rho_{jn}).$$
It follows that each superclass is in fact a single conjugacy class.  From this, the supercharacters are each irreducible.  These facts also follow from Corollary \ref{NoFourChains}, below.

The character formula of Theorem \ref{PatternGroupCharacterFormula} is simple for this example.  Vectors and matrices $M$, $a$, and $b$ can be indexed by $J$ in the linear order of (\ref{HeisenbergOrder}) above.  Note that for any $\phi, \eta\in J^*$ the matrix $M_\phi^\eta=0$ since there is never $(i,j,k,l)\in \cP$.   Note that
$$(a_\phi^\eta)_{ij}=\left\{\begin{array}{ll} \eta_{1n}\phi_{jn}, & \text{if $(i,j)=(1,j)$,}\\ 0, & \text{otherwise},\end{array}\right.\qquad \text{and}\qquad (b_\phi^\eta)_{jk}=\left\{\begin{array}{ll} \eta_{1n}\phi_{1j}, & \text{if $(j,k)=(j,n)$,}\\ 0, & \text{otherwise},\end{array}\right.$$
Now $\phi$ and $\eta$ mesh if $a_\phi^\eta=b_\phi^\eta=0$.   Thus, $\phi$ and $\eta$ mesh if and only if either
\begin{enumerate}
\item[(a)]  $\eta_{1n}=0$, or
\item[(b)]  $\phi_{ij}=0$ for all positions $(i,j)\neq (1,n)$.  
\end{enumerate}
In both cases, $\rank(M_\phi^\eta)=0$.  By Corollary \ref{CompleteLeftRightCoAction}, 
$$\corank(\eta)=\rank(M_L^\eta)=\left\{\begin{array}{ll} 0, & \text{if $\eta_{1n}=0$,}\\ n-2, & \text{otherwise.}\end{array}\right.$$
The vector $b$ in Theorem \ref{PatternGroupCharacterFormula} can be taken to be $0$ ($M_\phi^\eta b=-a_\phi^\eta$), so
\begin{equation}\label{HeisenbergCharacterFormula}
\overline{\chi^\eta(x_\phi)}=\left\{\begin{array}{ll} \prod_{(i,j)\in \cP} \theta(\phi_{ij}\eta_{ij}), & \text{if $\eta_{1n}=0$,}\\ q^{n-2} \theta(\phi_{1n}\eta_{1n}), & \text{if $\eta_{1n}\neq 0$, $x_\phi\in Z(H_n)$,}\\ 0, & \text{otherwise.}\end{array}\right.\end{equation}
Note that it is natural to choose orbit representatives for cosets of $Z(H_n)$ with $\eta_{1n}=0$.  Then (\ref{HeisenbergCharacterFormula}) gives the usual formula for the irreducible characters of $H_n$.

\subsubsection*{Example 2: $U_n(\FF_q)$.}  For the full group of upper-triangular matrices, $U_n(\FF_q)$, we may reduce the character formula to the one found in \cite{ADS04}. Let  $J=R^+$.  In this case, we may choose our orbit representatives $\phi\in J^*$ so that
\begin{equation}\label{PhiMonomialCondition} \phi_{jk}\neq 0\qquad\text{implies}\qquad \phi_{jl}=0=\phi_{ik}, \qquad\text{for all $j<l\neq k$ and $j\neq i<k$.}\end{equation}
Similarly, choose the co-orbit representatives $\eta\in J^*$ such that
\begin{equation}\label{EtaMonomialCondition} \eta_{jk}\neq 0\qquad\text{implies}\qquad \eta_{jl}=0=\eta_{ik}, \qquad\text{for all $j<l\neq k$ and $j\neq i<k$.}\end{equation}
That is, in either case we are permitted at most one nonzero entry for each row and column of $X_\phi, X_\eta\in \fkn$.  For such choices we may now compute $\chi^\eta(x_\phi)$.  

By our choice of representatives $\eta$ and $\phi$, the matrix $M_\phi^\eta$ has at most one nonzero entry in every row and column.  Suppose $(M_\phi^\eta)_{(ij)(kl)}\neq 0$ so that $(i,j,k,l)\in \cP$, $\phi_{jk}\neq 0$, and $\eta_{il}\neq 0$.  Since 
$$(a_\phi^\eta)_{(ij)}=\sum_{(i,j,k)} \phi_{jk}\eta_{ik},$$
our choice of representatives implies that in each summand either $\phi_{jk}=0$ or $\eta_{ik}=0$.  Thus, for every row of $M_\phi^\eta$ that has a nonzero row, the corresponding entry in $a_\phi^\eta$ is zero, and a solution to the matrix equation
$$M_\phi^\eta b_0=-a_\phi^\eta,$$
exists (and can be $b_0=0$) if and only if $a_\phi^\eta=0$.  Note that $a_\phi^\eta=0$ implies that if $(i,j,k)\in \cP$, then either $\phi_{jk}=0$ or $\eta_{ik}=0$.  In other words, any entry in $\supp(\eta)$ that shares a column with an entry of $\supp(\phi)$ must be below that entry. 

The nullspace $\Null(M_\phi^\eta)$ of $M_\phi^\eta$ has as a basis
$$\Null(M_\phi^\eta)=\CC\spanning\{e_{(k,l)}\ \mid\ \text{the $(k,l)$ column of $M_\phi^\eta$ is zero}\},$$
where $e_{(k,l)}$ is the vector with one in the $(k,l)$-slot and zeroes elsewhere.  Thus, for $b_\phi^\eta$ to be perpendicular to $\Null(M_\phi^\eta)$ it suffices that $(b_\phi^\eta)_{(k,l)}=0$ for every column $(k,l)$ of $M_\phi^\eta$ which has no nonzero entry.  Since
$$(b_\phi^\eta)_{k,l}=\sum_{(i,k,l)} \phi_{ik} \eta_{il},$$
this condition implies that for all $(i,k,l)\in \cP$, either $\phi_{ik}=0$ or $\eta_{il}=0$.  In other words, any entry in $\supp(\eta)$ that shares a row with an entry of $\supp(\phi)$ must be to the left of that entry. 

We may conclude $\chi^\eta(x_\phi)$ is nonzero ($\eta$ meshes with $\phi$) if and only if the following two conditions hold
\begin{enumerate}
\item[(a)] $\phi_{ij}\neq 0$ and $\eta_{il}\neq 0$ implies $j\geq l$,
\item[(b)] $\phi_{jk}\neq 0$ and $\eta_{ik}\neq 0$ implies $i \leq j$.
\end{enumerate}
This gives us a combinatorial interpretation of everything in the character formula except for the power of $q$.

By the choice of $\eta$, the matrix $M_R^\eta$ has at most one nonzero entry in every row and column.  Thus, the $\corank(\eta)=\rank(M_R^\eta)$ is the number of nonzero entries in $M_R^\eta$, which is 
$$|\{(i,j,k)\in \cP\ \mid\ \eta_{ik}\neq 0\}|=\sum_{(i,k)\in \cP\atop \eta_{ik}\neq 0} k-i-1.$$
In terms of matrices $k-i-1$ is the number of entries below the $(i,k)$ entry and above the diagonal. 

The rank of $M_\phi^\eta$ is also the number of nonzero entries in $M_\phi^\eta$, so
$$\rank(M_\phi^\eta)=|\{(i,j,k,l)\in \cP\ \mid\ \phi_{jk}\neq 0, \eta_{il}\neq 0\}|.$$
In terms of matrices, 
$$\corank(\eta)-\rank(M_\phi^\eta)=\hspace{-.25cm}\sum_{(i,l)\in \supp(\eta)} \left(\begin{array}{c} \text{Number of entries}\\ \text{below $(i,l)$ and}\\ \text{above the diagonal}\end{array}\right)-\left(\begin{array}{c} \text{Number of $(j,k)\in \supp(\phi)$}\\ \text{such that $(j,k)$ is}\\ \text{strictly SouthWest of $(i,l)$}\end{array}\right).$$ 

In the language of posets, if $\cP_{ij}$ is the the interval in $\cP$ from $i$ to $j$, then
$$\corank(\eta)-\rank(M_\phi^\eta)=\sum_{(i,l)\in \supp(\eta)} |\cP_{i+1,l-1}|-|\supp(\phi)\cap\cP_{i+1,l-1}|,$$
so
$$\overline{\chi^\eta(x_\phi)} =\left\{\begin{array}{ll} \hspace{-.15cm}\dd\prod_{(i,l)\in \supp(\eta)} \hspace{-.15cm} q^{|\cP_{i+1,l-1}|-|\supp(\phi)\cap\cP_{i+1,l-1}|}\theta(\eta_{il}\phi_{il}), &\text{if}\begin{array}{l}\text{$\phi_{ij}\neq 0$, $\eta_{il}\neq 0$ implies $j\geq l$,}\\ \text{$\phi_{jk}\neq 0$, $\eta_{ik}\neq 0$ implies $i \leq j$,}\end{array}\\ 0, & \text{otherwise}.\end{array}\right.$$

\subsection*{Example 3: No $4$-chains}

Let $\cP$ be any poset that contains no sequence of elements $i<j<k< l$, and for expository purposes assume that every element is contained in some $3$-chain.  Let
\begin{align*}
\mathcal{T} &=\{k\in \cP\ \mid\ (i,j,k)\in \cP\text{ for some $i,j\in \cP$}\}\\
\mathcal{M} &=\{j\in \cP\ \mid\ (i,j,k)\in \cP\text{ for some $i,k\in \cP$}\}\\
\mathcal{B} &= \{i\in \cP\ \mid\ (i,j,k)\in \cP\text{ for some $j,k\in \cP$}\}
\end{align*}
so that by our assumption $\cP=\cT\cup\cM\cup\cB$.
The group $U_\cP$ consists of matrices of the form
$$u(A,C,B)=\left(\begin{array}{ccc} Id_{|\cB|} & A & C \\ 0 & Id_{|\cM|} &  B\\ 0 & 0 & Id_{|\cT|}\end{array}\right).$$
Note that if $\cB<\cM<\cT$ in $\cP$, then $U_\cP$ is the unipotent radical of the parabolic subgroup 
$$P_{|\cB|,|\cM|,|\cT|}=\left\{\left(\begin{array}{ccc} A &  \ast  & \ast\\ 0 & B & \ast \\ 0 & 0 & C\end{array}\right)\ \mid\ A\in \GL_{|\cB|}(\FF_q), B\in \GL_{|\cM|}(\FF_q),C\in \GL_{|\cT|}(\FF_q)\right\}$$
of $\GL_{|\cB|+|\cM|+|\cT|}(\FF_q)$.

If $u(A,C,B)\in U_\cP$, then the corresponding two-sided orbit is 
$$\left\{\left(\begin{array}{c@{}c@{}c} Id_{|\cB|} & A & C+AX+YB \\ 0 & Id_{|\cM|} &  B\\ 0 & 0 & Id_{|\cT|}\end{array}\right)
\bigg|
\left(\begin{array}{c@{}c@{}c} Id_{|\cB|} & 0 & 0 \\ 0 & Id_{|\cM|} &  X\\ 0 & 0 & Id_{|\cT|}\end{array}\right), \left(\begin{array}{c@{}c@{}c} Id_{|\cB|} & Y & 0 \\ 0 & Id_{|\cM|} &  0\\ 0 & 0 & Id_{|\cT|}\end{array}\right)\in U_\cP\right\}.$$
Similarly, if $\lambda(U,W,V)\in \fkn_\cP^*$ corresponds to the matrix $u(U,W,V)$, then the co-orbit containing $\lambda(U,W,V)$ is
$$\{\lambda(U+W(X^{\Tr}), W, V+(Y^{\Tr}) W)\ \mid\ \left(\begin{array}{c@{}c@{}c} Id_{|\cB|} & 0 & 0 \\ 0 & Id_{|\cM|} &  -X\\ 0 & 0 & Id_{|\cT|}\end{array}\right), \left(\begin{array}{c@{}c@{}c} Id_{|\cB|} & -Y & 0 \\ 0 & Id_{|\cM|} &  0\\ 0 & 0 & Id_{|\cT|}\end{array}\right)\in U_\cP\}.$$

Since there are no sequences of length four, the matrix $M_\phi^\eta=0$, so $\phi$ meshes with $\eta$ if and only if $a_\phi^\eta=b_\phi^\eta=0$.   But,
$$(a_\phi^\eta)_{ij}=\left\{\begin{array}{ll} \dd\sum_{(j,k)\in \cP} \phi_{jk}\eta_{ik}, & \text{$i\in \cB,j\in \cM$,}\\ 0, & \text{otherwise,}\end{array}\right.
\quad\text{and}\quad
(b_\phi^\eta)_{jk}=\left\{\begin{array}{ll} \dd\sum_{(i,j)\in\cP} \phi_{ij}\eta_{ik}, & \text{$j\in \cM,k\in \cT$,}\\ 0, & \text{otherwise,}\end{array}\right.$$
which translate into the conditions for $u(A,C,B)$ and $\lambda(U,W,V)$, 
\begin{equation}\label{NoFourChainsMeshing}
A_{\ast j}\cdot W_{\ast k}=0, \quad \text{if $j<k$},\quad\text{and}\quad B_{j\ast}\cdot W_{i\ast}=0, \quad\text{if $i<j$,}
\end{equation}
where $M_{\ast j}$ is the $j$th column of $M$, $M_{i\ast}$ is the $i$th row of $M$ and $\cdot$ is the usual vector dot product.  For $j\in \cM$, let
$$W^{(j)}=\text{columns $W_{\ast k}$ of $W$ that $j<k$ in $\cP$}.$$ 

The resulting character formula is
\begin{equation}\label{No4ChainFormula}\chi^{\lambda(U,W,V)}(u(A,C,B))=\left\{\begin{array}{ll} \dd q^{\sum_{j\in \cM} \rank(W^{(j)})}\prod_{(i,j)\in \cP} \theta(\lambda_{ij} u_{ij}), & \text{if $\lambda$ and $u$ satisfy (\ref{NoFourChainsMeshing}),}\\ 0, & \text{otherwise}.\end{array}\right.\end{equation}
Proposition \ref{SuperNormal} below shows that for these examples supercharacters and superclasses are irreducible characters and conjugacy classes so that (\ref{No4ChainFormula}) gives a formula for the irreducible characters of $U_\cP$.

\subsection{Supertheory versus usual character theory}

For $\eta\in J^*$, let
\begin{align*}
\ann_J^R(\eta)&=\{\rho\in J^*\ \mid\ X_\phi X_\rho\in \ker(\lambda_\eta),\text{ for all $\phi\in J^*$}\}\\
\ann_J^L(\eta)&=\{\rho\in J^*\ \mid\ X_\rho X_\phi \in \ker(\lambda_\eta),\text{ for all $\phi\in J^*$}\}
\end{align*}
Diaconis and Isaacs give the following characterization of which supercharacters are irreducible characters \cite{DI06}.

\begin{theorem}[Diaconis, Isaacs] \label{DICharacters}
Let $J$ be closed, and $\eta\in J^*$.  Then $\chi^\eta$ is irreducible if and only if $\ann_J^R(\eta)+\ann_J^L(\eta)=J^*$.
\end{theorem}

This implies the following more combinatorial proposition.

\begin{proposition} \label{SuperNormal} Let $J\subseteq R^+$ be a closed subset.
\begin{enumerate}
\item[(a)] Suppose $\phi\in J^*$ and there is no $4$-chain $(i,j,k,l)\in \cP_J$ such that $(i,j),(k,l)\in \supp(\phi)$.  Then the two sided orbit of $x_\phi$ is a conjugacy class.
\item[(b)] Suppose $\eta\in J^*$ and there is no $4$-chain $(i,j,k,l)\in \cP_J$ such that $(i,k),(j,l)\in \supp(\eta)$.  Then $\chi^\eta$ is an irreducible character.
\end{enumerate}
\end{proposition}
\begin{proof}
(a) Note that if $\phi$ satisfies the hypothesis, then any $\phi'\in O_\phi$ satisfies the hypothesis.
 If the hypothesis is satisfied, then for every $(i,j)\in J$, the group element $x_{ij}(t)$ acts trivially either from the right or from the left.  Thus, conjugation by $x_{ij}(t)$ is either a right or left action.  

(b)  Suppose there is no $(i,j,k,l)\in \cP$ such that $(i,k),(j,l)\in \supp(\eta)$, and let $\rho\in J^*$.  By Theorem \ref{DICharacters}, it suffices to find $\rho'\in \ann_J^R(\eta)$ and $\rho''\in \ann_J^L(\eta)$ such that $\rho=\rho'+\rho''$.  Let
\begin{align*}
\rho'_{jk} &=\left\{\begin{array}{ll} 0, & \text{if there exists $(i,j,k)\in \cP$ with $(i,k)\in \supp(\eta)$,}\\ \rho_{jk}, &\text{otherwise.}\end{array}\right.\\
\rho_{jk}'' & = \left\{\begin{array}{ll} 0, & \text{if there exists $(j,k,l)\in \cP$ with $(j,l)\in \supp(\eta)$ OR $\rho_{jk}'\neq 0$,}\\ \rho_{jk}, &\text{otherwise.}\end{array}\right.
\end{align*}
Note that by the assumption, if $\rho_{jk}\neq 0$, then either $\rho_{jk}'\neq 0$ or $\rho_{jk}''\neq 0$.  Thus, $\rho=\rho'+\rho''$.

For $\phi\in J^*$,
$$X_\phi X_{\rho'} = -X_\phi + X_\phi x_{\rho'}.$$
If $X_\phi X_{\rho'}=X_{\phi'}$, then by Theorem \ref{LeftRightAction},
$$\phi_{il}'=\sum_{(i,j,l)\in \cP} \phi_{ij}\rho_{jl}',$$
so
$$\lambda^\eta(X_\phi X_{\rho'})=\sum_{(i,k)}\eta_{ik}\sum_{(i,j,k)}\phi_{ij}\rho_{jk}'=\sum_{(i,j)}\phi_{ij}\biggl(\sum_{(i,j,k)} \eta_{ik}\rho_{jk}'\biggr)=0.$$
Thus, $\rho'\in \ann_J^R(\eta)$.  Similarly, $\rho''\in \ann_J^L(\eta)$, so by Theorem \ref{DICharacters} $\chi^\eta$ is an irreducible character.
\end{proof}

\begin{corollary} \label{NoFourChains}
If $\cP_J$ has no $4$-chains, then all superclasses are conjugacy classes, and all supercharacters are irreducible characters.
\end{corollary}

The following Corollary uses \cite{DI06} and the fact that we may choose our co-orbit representatives to satisfy (\ref{EtaMonomialCondition}) to obtain a full characterization of when supercharacters are irreducible.

\begin{corollary} \label{FullUCharacters}
Let $J=R^+$ so that $U_J=U_n(\FF_q)$.  Let $\eta\in J^*$ such that $\eta$ satisfies (\ref{EtaMonomialCondition}).  Then  $\chi^\eta$ is an irreducible character if and only if  there is no $4$-chain $(i,j,k,l)\in \cP_J$ such that $(i,k),(j,l)\in \supp(\eta)$. 
\end{corollary}
\begin{proof}
Consider
\begin{align*}
\ann_J^R(\eta)&=\{\rho\in J^*\ \mid\ \sum_{(i,k)\in \supp(\eta)} \eta_{ik}\biggl(\sum_{(i,j,k)} \phi_{ij}\rho_{jk}\biggr)=0, \text{for all $\phi\in J^*$}\} \\
&=\{\rho\in J^*\ \mid\ \sum_{k\geq j}\eta_{ik}\rho_{jk}=0, \text{for all $(i,j)\in J$}\}\\
\ann_J^L(\eta)&=\{\rho\in J^*\ \mid\ \sum_{(i,k)\in \supp(\eta)} \eta_{ik}\biggl(\sum_{(i,j,k)} \rho_{ij}\phi_{jk}\biggr)=0, \text{for all $\phi\in J^*$}\}\\
&=\{\rho\in J^*\ \mid\ \sum_{i\leq j}\eta_{ik}\rho_{ij}=0, \text{for all $(j,k)\in J$}\}. 
\end{align*}
By assumption, $(i,j),(k,l)\in \supp(\eta)$ implies that $i\neq k$ and $j\neq l$, so
\begin{align*}
\ann_J^R(\eta)&=\{\rho\in J^*\ \mid\ \eta_{ik}\rho_{jk}=0, \text{for all $(i,j,k)\in \cP$, $(i,k)\in \supp(\eta)$}\}\\
\ann_J^L(\eta)&=\{\rho\in J^*\ \mid\ \eta_{ik}\rho_{ij}=0, \text{for all $(i,j,k)\in \cP$, $(i,k)\in \supp(\eta)$}\}. 
\end{align*}
Thus, 
$$\ann_J^R(\eta)+\ann_J^L(\eta)=J^*$$
if and only if there does not exists $(i,j,k)\in \cP$ and $(j,k,l)\in \cP$ such that $(i,k), (j,l)\in \supp(\eta)$.
\end{proof}

\noindent\textbf{Examples.} 
\begin{enumerate}
\item Corollary \ref{NoFourChains} implies that the supertheory for the Heisenberg group and the pattern groups of Section \ref{SuperCharacterFormula}, Example 3, is just the ordinary irreducible character theory.
\item  If $\phi, \eta\in J^*$ satisfy $|\supp(\phi)|=|\supp(\eta)|=1$, then $x_\phi$ is a conjugacy class and $\chi^\eta$ is an irreducible character. 
\item If $U_\cP=U_n(\FF_q)$, then the supercharacters $\chi^\eta$ with maximal degree are  given by  $\supp(\eta)=\{(i,n+1-i)\ \mid\ 1\leq i\leq \lfloor n/2\rfloor\}$, or  in terms of the poset, $\supp(\eta)$ is the set of pairs
$$
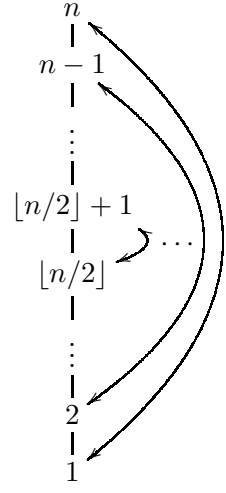
\xymatrix@R=.3cm@C=.25cm{ {n} \\ {n-1} \ar @{-} [u] \\ {\vdots} \ar @{-} [u]\\ {\lfloor n/2\rfloor+1} \ar @{-} [u] \\ {\lfloor n/2\rfloor} \ar @{-} [u] \ar @{<->} @/_{1cm}/ [u] \save []+<1.4cm,.4cm> *{\cdots}\restore \\ {\vdots} \ar @{-} [u] \\ {2} \ar @{-} [u] \ar @{<->} @/_{1.75cm}/ [uuuuu]  \\{1} \ar @{-} [u] \ar @{<->} @/_{2cm}/ [uuuuuuu]}$$
While the poset contains $i<j<k<l$, there is no such sequence such that $(i,k),(j,l)\in \supp(\eta)$.  Thus, Proposition \ref{SuperNormal} implies $\chi^\eta$ is an irreducible character.   Note that since the degrees of supercharacters are bigger than the degrees of irreducible characters, we have that the all irreducible characters of maximal degree are supercharacters of this form. 

These characters were first identified by Lehrer \cite{Le74}.  Andr\'e \cite{An02} and Yan \cite{Ya01} realized that they are in fact supercharacters.  The fact that they constitute all of the characters of maximal degree is first proved by Isaacs \cite{Is06}.

Note that if $x_\phi$ is such that $\supp(\phi)=\{(i,n+1-i)\ \mid\ 1\leq i\leq \lfloor n/2\rfloor\}$, then $x_\phi$ is also a conjugacy class.  
\end{enumerate}

\vspace{.25cm}

\noindent\textbf{(Counter) Examples.}  The following examples show that while the conditions in Proposition \ref{SuperNormal} are necessary, they are not sufficient.  

\vspace{.25cm}

\noindent 1. One can show that if
$$\cP=\ 
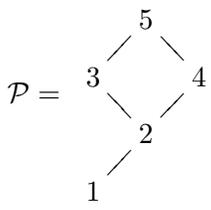
\xy<0cm,1cm>\xymatrix@R=.3cm@C=.3cm{*={} & {5} \\ {3} \ar @{-} [ur] & *={} & {4} \ar @{-} [ul]\\ *={} & {2}\ar @{-} [ul] \ar @{-} [ur] \\ {1} \ar @{-} [ur]}\endxy$$
and $\supp(\phi)=\{(1,2),(2,4),(3,5)\}$, then $O_\phi$ is a conjugacy class (even though $(1,2,3,5)\in \cP$).  See \cite{An01} for a necessary and sufficient condition in the case when $\cP$ is a chain.

\vspace{.25cm}

\noindent 2. Suppose
$$\cP=
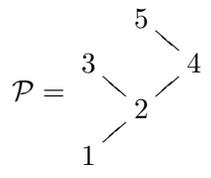
\xy<0cm,1cm>\xymatrix@R=.15cm@C=.3cm{*={} & {5} \\   {3} & *={} & {4}\ar @{-} [ul] \\ *={} & {2} \ar @{-} [ur] \ar @{-} [ul] \\   {1} \ar @{-} [ur] }\endxy$$
and $\supp(\eta)=\{(1,3),(1,4),(2,5)\}$.  Then, 
\begin{align*}
 \ann_J^R(\eta) &= \{\rho'\in J^*\ \mid\ \eta_{13}\phi_{12}\rho_{23}'+\eta_{14}\phi_{12}\rho_{24}'+\eta_{25}\phi_{24}\rho_{45}'=0, \text{ for all $\phi_{12},\phi_{24}\in \FF_q$}\}\\
 &=  \{\rho'\in J^*\ \mid\ \rho_{45}'=0, \rho_{23}'=-\eta_{13}^{-1}\rho_{24}'\eta_{14}\},\\
 \ann_J^R(\eta) &= \{\rho''\in J^*\ \mid\ \eta_{13}\rho_{12}''\phi_{23}+\eta_{14}\rho_{12}''\phi_{24}+\eta_{25}\rho_{24}''\phi_{45}=0, \text{ for all $\phi_{23},\phi_{24}, \phi_{45}\in \FF_q$}\}\\
 &=  \{\rho''\in J^*\ \mid\ \rho_{12}''=\rho_{24}''=0\}.
\end{align*}
Given $\rho\in J^*$, define $\rho'\in \ann_J^R(\eta)$ and $\rho''\in \ann_J^L(\eta)$ by
\begin{align*} \rho_{ij}' &=\left\{\begin{array}{ll} \rho_{ij}, & \text{if $(i,j)\notin \{(2,3),(4,5)\}$,}\\ 0, & \text{if $(i,j)=(4,5)$,}\\ -\eta_{13}^{-1}\rho_{24}\eta_{14}, & \text{if $(i,j)=(2,3)$.}\end{array}\right.\\
\rho_{ij}'' &=\left\{\begin{array}{ll} 0, & \text{if $(i,j)\notin \{(2,3),(4,5)\}$,}\\ \rho_{45}, & \text{if $(i,j)=(4,5)$,}\\ \rho_{23}+\eta_{13}^{-1}\rho_{24}\eta_{14}, & \text{if $(i,j)=(2,3)$.}\end{array}\right.
\end{align*}
Since $\rho=\rho'+\rho''$, Theorem \ref{DICharacters} implies that $\chi^\eta$ is an irreducible character. 

\vspace{.25cm}

\noindent 3. This example shows that determining whether a supercharacter $\chi^\eta$ is irreducible can require more information than just the support of $\eta$.  Suppose
$$\cP=
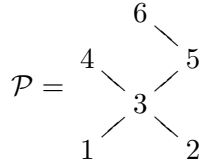
\xy<0cm,1cm>\xymatrix@R=.15cm@C=.3cm{*={} & {6} & *={} \\ {4} &   *={}   & {5}  \ar @{-} [ul]  \\ *={} &  {3}  \ar @{-} [ur]   \ar @{-} [ul] \\   {1} \ar @{-} [ur]  & *={} & {2} \ar @{-} [ul] }\endxy$$
and $\supp(\eta)=\{(1,4),(1,5),(2,4),(2,5),(3,6)\}$. Then,
\begin{align*}
\hspace*{-.25cm} \ann_J^R(\eta) &= \{\rho\in J^*\ \mid\ \eta_{14}\phi_{13}\rho_{34}+\eta_{15}\phi_{13}\rho_{35}+\eta_{24}\phi_{23}\rho_{34}+\eta_{25}\phi_{23}\rho_{35}+\eta_{36}\phi_{35}\rho_{56}=0, \phi\in J^*\}\\
 &=\{\rho\in J^*\ \mid\ \rho_{56}=0, \eta_{14}\rho_{34}=-\eta_{15}\rho_{35}, \eta_{24}\rho_{34}=-\eta_{25}\rho_{35}\},\\
 \ann_J^R(\eta) &= \{\rho\in J^*\ \mid\ \eta_{14}\rho_{13}\phi_{34}+\eta_{15}\rho_{13}\phi_{35}+\eta_{24}\rho_{23}\phi_{34}+\eta_{25}\rho_{23}\phi_{34}+\eta_{36}\rho_{35}\phi_{56}=0,\phi\in J^*\}\\
 &=  \{\rho\in J^*\ \mid\ \rho_{35}=0, \eta_{14}\rho_{13}=-\eta_{24}\rho_{23}, \eta_{15}\rho_{13}=-\eta_{25}\rho_{23} \}.
\end{align*}
Let $\rho\in J^*$.  We want $\rho=\rho'+\rho''$ for $\rho'\in  \ann_J^R(\eta)$ and $\rho''\in  \ann_J^L(\eta)$.  Since elements of $\ann_J^R(\eta)$  can have arbitrary values at $J\setminus\{(5,6),(3,4), (3,5)\}$, and elements of  $\ann_J^L(\eta)$ can have arbitrary values at $(5,6)$ and $(3,4)$, we can find $\rho'$ and $\rho''$ such that $\rho_{ij}=\rho_{ij}'+\rho_{ij}''$ for $(i,j)\neq (3,5)$.  However, $\rho_{35}''=0$, so 
$$\rho_{35}=\rho_{35}'=-\eta_{15}^{-1}\eta_{14}\rho_{34}'=-\eta_{25}^{-1}\eta_{24}\rho_{34}',$$
which can happen if and only if $\eta_{25}\eta_{14}-\eta_{15}\eta_{24}=0$.  Thus, $\chi^\eta$ is an irreducible character if and only if $\eta_{25}\eta_{14}-\eta_{15}\eta_{24}=0$.

\section{An algebra group supercharacter formula}

Let $H\subseteq U_n(\FF_q)$ be an algebra group with corresponding vector space $V_H$ (as in Proposition \ref{AlgebraGroupCharacterization}), and let $\fkn_H=H-1$.  This section gives a formula for a general supercharacter on a general superclass.  The main result is followed by a corollary regarding supercharacter values and by several examples.

Fix a basis $\{v_1,v_2,\ldots,v_d\}$ for $V_H$ and let $\{\lambda_1,\lambda_2,\ldots,\lambda_d\}$ be a dual basis for the dual space $V_H^*$.   For $\phi=(\phi_1,\phi_2,\ldots, \phi_d)\in \FF_q^d$, let $x_\phi\in H$ be the element that corresponds to the vector
$$\phi_1v_1+\phi_2v_2+\cdots+\phi_d v_d\in V_H.$$
Similarly, for $\eta=(\eta_1,\eta_2,\ldots,\eta_d)\in \FF_q^d$, let $\lambda_\eta\in \fkn_H^*$ be the functional corresponding to 
$$\eta_1\lambda_1+\eta_2\lambda_2+\cdots+\eta_d\lambda_d\in V_H^*.$$

Let $c_{ij}^k\in \FF_q$ be defined by
$$X_{v_i}X_{v_j}=\sum_{k=1}^d c_{ij}^k X_{v_k},$$
where $1+X_{v_i}\in H$ is the group element corresponding to the basis vector $v_i\in V_H$.
Let $C_i$ and $C^j$ denote the $d\times d$ matrices given by
$$(C_i)_{jk}=c_{ij}^k \qquad\text{and}\qquad (C^j)_{ik}=c_{ij}^k.$$
For $\phi,\eta\in \FF_q^d$, let $M_\phi^\eta$ be the $d\times d$ matrix given by
\begin{equation} (M_\phi^\eta)_{ij}=\phi C_i C^j\eta. \label{AlgebraGroupMatrix}\end{equation}
Define $a_\phi^\eta,b_\phi^\eta\in \FF_q^d$, by 
\begin{align}
(a_\phi^\eta)_i &=\phi C_i\eta \label{AlgebraGroupaVector}\\
(b_\phi^\eta)_j &= \phi C^j \eta\label{AlgebraGroupbVector}
\end{align}

We say that $\phi$ \emph{meshes} with $\eta$ if 
\begin{enumerate}
\item[(1)] there exists $b\in \FF_q^d$ such that $M_\phi^\eta b=-a_\phi^\eta$,
\item[(2)] $b_\phi^\eta$ is perpendicular to $\Null(M_\phi^\eta)$.
\end{enumerate}

\begin{theorem} \label{AlgebraGroupCharacterFormula} Let $H$ be an algebra group with $\dim(V_H)=d$.  
Let $\phi,\eta\in \FF_q^d$ so that $x_\phi\in H$ and $\lambda_\eta\in \fkn_H^*$ are defined according to bases of $V_H$ and $V_H^*$, respectively.  Let $M_\phi^\eta$, $a_\phi^\eta$, and $b_\phi^\eta$ be as in (\ref{AlgebraGroupMatrix}), (\ref{AlgebraGroupaVector}), and (\ref{AlgebraGroupbVector}), respectively.  Then
$$\overline{\chi^\eta(x_\phi)}= \left\{\begin{array}{ll} \dd q^{\mathrm{corank}(\eta)-\mathrm{rank}(M_\phi^\eta)} \theta(b_0\cdot b_\phi^\eta)  \prod_{i=1}^d  \theta(\phi_i\eta_i),  & \text{if  $\phi$ meshes with $\eta$,}\\ 0, & \text{otherwise,}\end{array}\right.$$
where $b_0\in \FF_q^d$ satisfies $M_\phi^\eta b_0=-a_\phi^\eta$, and $q^{\corank(\eta)}$ is the size of the right $H$-orbit of $\lambda_\eta$.
\end{theorem}

\begin{proof}
From the definitions in Section \ref{SectionIntroSuperCharacters},
$$\overline{\chi^\eta(x_\phi)} = \frac{q^{\corank(\eta)}}{|O_\phi|}\sum_{\rho\in \FF_q^d} \theta\circ\lambda_\eta(X_\rho)=\frac{q^{\corank(\eta)}}{q^{2d}}\sum_{a,b\in \FF_q^d} \theta\circ\lambda_\eta(x_a X_\phi x_b)$$
Note that for $a=(a_1,a_2,\ldots, a_d),b=(b_1,b_2,\ldots,b_d)\in \FF_q^d$,
$$x_a X_\phi x_b=\sum_{i=1}^d \phi_i X_{v_i} +\sum_{i,j,k} (a_i\phi_j+\phi_i b_j) c_{ij}^k X_{v_k} + \sum_{i,j,k,l,m} a_i\phi_j b_l c_{ij}^k c_{kl}^m X_{v_m}.$$
Thus,
\begin{align*}
\overline{\chi^\eta(x_\phi)} &= \frac{q^{\corank(\eta)}}{q^{2d}}\sum_{a,b\in \FF_q^d} \theta\left(\sum_{m=1}^d \eta_m\left(\phi_m+\sum_{i,j} (a_i\phi_j+\phi_i b_j) c_{ij}^m + \sum_{i,j,k,l} a_i\phi_j b_l c_{ij}^k c_{kl}^m \right)\right)\\
&=\frac{q^{\corank(\eta)}\theta_{\phi\eta}}{q^{2d}}\sum_{a,b\in \FF_q^d} \theta\left(\sum_{m=1}^d \eta_m\left(\sum_{i,j} (a_i\phi_j+\phi_i b_j) c_{ij}^m + \sum_{i,j,k,l} a_i\phi_j b_l c_{ij}^k c_{kl}^m \right)\right),
\end{align*}
where
$$\theta_{\phi\eta}=\theta(\phi_1\eta_1+\phi_2\eta_2+\cdots+\phi_d\eta_d).$$
Collect all terms involving $a_i$ to obtain
\begin{align*}
\overline{\chi^\eta(x_\phi)} &=\frac{q^{\corank(\eta)}\theta_{\phi\eta}}{q^{2d}}\sum_{a,b\in \FF_q^d}\theta\left( \sum_{i=1}^d a_i\biggl(\sum_{j,m}\phi_jc_{ij}^m \eta_m+\sum_{j,k,l,m} \phi_jc_{ij}^k c_{kl}^m\eta_m b_l\biggr) +\sum_{i,j,m} \phi_i c_{ij}^m \eta_m b_j \right)\\
&= \frac{q^{\corank(\eta)}\theta_{\phi\eta}}{q^{2d}}\sum_{a,b\in \FF_q^d}\theta\left( a\cdot(a_\phi^\eta+M_\phi^\eta b) +b_\phi^\eta \cdot b \right).
\end{align*}
Note that 
$$\sum_{a\in \FF_q^d} \theta( a\cdot(a_\phi^\eta+M_\phi^\eta b))=0,$$
unless $M_\phi^\eta b=-a_\phi^\eta$.  Thus, if $\mathcal{S}=\{b\in \FF_q^d\ \mid\ M_\phi^\eta b=-a_\phi^\eta\}$, then
$$\overline{\chi^\eta(x_\phi)} = \frac{q^{\corank(\eta)}\theta_{\phi\eta}}{q^{d}}\sum_{b\in \mathcal{S}}\theta\left(b_\phi^\eta \cdot b \right).$$
Fix an element $b_0\in \mathcal{S}$.  Since every other vector in $\mathcal{S}$ is of the form $b_0+b'$, where $b'$ is in the nullspace of $M_\phi^\eta$, we have
$$\overline{\chi^\eta(x_\phi)} = \frac{q^{\mathrm{corank}(\eta)}\theta_{\phi\eta}\theta(b_0\cdot b_\phi^\eta)}{q^{d}}\sum_{b'\in \mathrm{Null}(M_\phi^\eta)}\theta(b'\cdot b_\phi^\eta).$$
Let $\{b_1',b_2',\ldots, b_r'\}$ be a basis for $\mathrm{Null}(M_\phi^\eta)$.  Then
\begin{equation*}
\sum_{b'\in \mathrm{Null}(M_\phi^\eta)}\theta(b'\cdot b_\phi^\eta)=\sum_{t\in \FF_q^r} \prod_{i=1}^r \theta(t_i b_i'\cdot b_\phi^\eta)= \prod_{i=1}^r \sum_{t_i\in \FF_q}\theta(t_i b_i'\cdot b_\phi^\eta).
\end{equation*}
Thus, if any $b_i'$ is not orthogonal to $b_\phi^\eta$, then the product is zero.  If $\phi$ meshes with $\eta$, then
$$\overline{\chi^\eta(x_\phi)} = \frac{q^{\mathrm{corank}(\eta)}\theta_{\phi\eta}\theta(b_0\cdot b_\phi^\eta)}{q^{d}}q^{d-\rank(M_\phi^\eta)}=q^{\corank(\eta)-\rank(M_\phi^\eta)}\theta_{\phi\eta}\theta(b_0\cdot b_\phi^\eta),$$
as desired.
\end{proof}

\begin{corollary}\label{SuperCharacterValues}
Let $H\subseteq U_n(\FF_q)$ be an algebra group and suppose $\text{char}(\FF_q)=p$.  Then
\begin{enumerate}
\item[(a)]  The nonzero supercharacter values are integer multiples of $p$th roots of unity,
\item[(b)]  If $p=2$, then all supercharacters are real-valued.
\end{enumerate}
\end{corollary}

\noindent\textbf{Remarks.}  
\begin{enumerate}
\item Note that if $H=U_J$ is a pattern group, then a natural basis of $V_H$ is $\{v^{(ij)}\ \mid\ (i,j)\in J\}$, where $v^{(ij)}\in \FF_q^{|J|}$ is given by $v^{(ij)}_{rs}=\delta_{ij,rs}$.  With this basis, Theorem \ref{AlgebraGroupCharacterFormula} reduces to Theorem \ref{PatternGroupCharacterFormula}.  
\item Corollary \ref{SuperCharacterValues} (a) is the supercharacter version of a conjecture by Isaacs \cite[Conjecture 10.1]{Is06} concerning irreducible characters.
\end{enumerate}

\subsection{Examples}  

\subsubsection*{1. A supercharacter formula}

Let $H_{n-1}=1+X\FF_q[X]/(X^{n-1})$ be the abelian group of polynomials that begin with 1 under polynomial multiplication mod $X^{n-1}$.  The supercharacters for this algebra group were computed as an example in \cite[Appendix C]{DI06}.  We extend these calculations by considering the semi-direct product of $H_{n-1}$ with $\FF_q^{n-1}$, the  abelian group under vector addition.  

For $f=1+\sum_{i=2}^{n-1} a_{i} X^{i-1}\in H_{n-1}$, $t=(t_1,t_2,\ldots, t_{n-1})\in \FF_q^{n-1}$, define $f t f^{-1}=(t'_1,t_2',\ldots,t_{n-2}')\in \FF_q^{n-1}$ by 
$$t_j'=t_j+\sum_{k=j+1}^{n-1} a_k t_k.$$
This action makes $H=H_{n-1}\ltimes \FF_q^{n-1}$ an algebra group.  In fact, as a subgroup of $U_n$,
$$H=\left\{\left(
\begin{array}{ccccc|c} 
1 & a_2 & a_3 & \cdots & a_{n-1} & t_1\\ 
0 & \ddots & \ddots & \ddots & \vdots & \vdots\\
 &  & 1 & a_2 & a_3 & t_{n-3}\\
\vdots &  & \ddots & 1 & a_2 & t_{n-2}\\
0 & 0 & \cdots & 0 & 1 & t_{n-1}\\ \hline
0 & 0 & \cdots & 0 & 0 & 1\end{array}\right)\ \bigg| a_2,\ldots, a_{n-1},t_1,\ldots, t_{n-1}\in \FF_q\right\}.$$ 
Let $\{v_2,v_3,\ldots,v_{n-1},v_{\bar{1}},v_{\bar{2}},\ldots,v_{\overline{n-1}}\}$ be the basis of $V_H$ given by 
\begin{align*}
X_{v_i}&=X\in \fkn, & & \text{where $X_{1i}=X_{2,i+1}=\cdots = X_{n-i,n-1}=1$ with zeroes elsewhere,}\\
X_{v_{\bar{j}}} & = X\in \fkn, & & \text{where $X_{jn}=1$ with  zeroes elsewhere.} 
\end{align*}
Let $\{\lambda_2,\ldots, \lambda_{n-1},\lambda_{\bar{1}},\ldots,\lambda_{\overline{n-1}}\}$ be the corresponding dual basis of $V_H^*$.

A small computation shows that the superclasses are indexed by 
$$\{(i,\bar{j},a,t),(i,\cdot,a,0),(\cdot,\bar{j},0,t),(\cdot,\cdot,0,0)\ \mid\ 2\leq i\leq n-1, n-i< j\leq n-1, a,t\in \FF_q^\times\},$$
such that natural representatives for the orbits $O_{(i,\bar{j},a,t)}$, $O_{(i,\cdot,a,0)}$, $O_{(\cdot,\bar{j},0,t)}$,and $O_{(\cdot,\cdot,0,0)}$ correspond to the vectors
$$av_i+bv_{\bar{j}},\quad av_i,\quad bv_{\bar{j}},\quad 0,\qquad \text{respectively}.$$
Similarly, the co-orbits are indexed by the same quadruplets such that a natural representative for the orbits $O^{(i,\bar{j},a,t)}$, $O^{(i,\cdot,a,0)}$, $O^{(\cdot,\bar{j},0,t)}$,and $O^{(\cdot,\cdot,0,0)}$ are the functionals corresponding to the vectors
$$a\lambda_i+b\lambda_{\bar{j}},\quad a\lambda_i, \quad b\lambda_{\bar{j}},\quad 0,\qquad \text{respectively}.$$

We first calculate the basic ingredients $M_\phi^\eta$, $a_\phi^\eta$, and $b_\phi^\eta$ of (\ref{AlgebraGroupMatrix}), (\ref{AlgebraGroupaVector}), and (\ref{AlgebraGroupbVector}).  Note that 
\begin{align*} 
c_{ij}^k& = \left\{\begin{array}{ll} 1, &\text{if $i+j\leq k$, $k\in \{2,\ldots, n-1\}$,}\\ 0, & \text{otherwise,}\end{array}\right. & & c_{\bar{i}\bar{j}}^k =0,\\
c_{i\bar{j}}^k & =\left\{\begin{array}{ll} 1, &\text{if $k=j-i+1$, $k\in\{\bar{1},\ldots,\overline{n-1}\}$,}\\ 0, & \text{otherwise,}\end{array}\right. & & c_{\bar{i}j}^k  =0,
\end{align*} 
so that by the definition (\ref{AlgebraGroupMatrix}), for $a,b,s,t\in \FF_q$, the $f,g$ entry of $M_{(i,\bar{j},a,s)}^{(k,\bar{l},b,t)}$ is
\begin{align*}
&(ac_{fi}^{2}+sc_{f\bar{j}}^2,\ldots,ac_{fi}^{n-1}+sc_{f\bar{j}}^{n-1},ac_{fi}^{\bar{1}}+sc_{f\bar{j}}^{\bar{1}},\ldots, ac_{fi}^{\overline{n-1}}+sc_{f\bar{j}}^{\overline{n-1}})\left(\begin{array}{@{}c@{}} b c_{2g}^{k}+t c_{2g}^{\bar{l}}\\ 
\vdots\\
b c_{n-1,g}^{k}+t c_{n-1,g}^{\bar{l}} \\
 0\\
 \vdots\\
0 \end{array}\right)\\
 &=0, \qquad \text{unless $f\in  \{2,\ldots, n-1\}$}.
 \end{align*}
If $f,g\in \{2,\ldots, n-1\}$, then
\begin{align*} (M_{(i,\bar{j},a,s)}^{(k,\bar{l},b,t)})_{f,g} & = (0,\ldots, 0, ac_{fi}^{f+i},0,\ldots, 0, sc_{f\bar{j}}^{\overline{j-f+1}},0,\ldots,0)\left(\begin{array}{c} 0\\ \vdots \\ 0\\ bc_{k-g,g}^k\\ 0\\ \vdots\\ 0 \end{array}\right)\\
&= \left\{\begin{array}{ll} ab, & \text{if $k-g=f+i$,}\\ 0, & \text{otherwise}.\end{array}\right.
\end{align*}
If $f\in \{2,\ldots, n-1\}$ and $g\in \{\bar{1},\ldots, \overline{n-1}\}$, then
\begin{align*} (M_{(i,\bar{j},a,s)}^{(k,\bar{l},b,t)})_{f,g} & = (0,\ldots, 0, ac_{fi}^{f+i},0,\ldots, 0, sc_{f\bar{j}}^{\overline{j-f+1}},0,\ldots,0)\left(\begin{array}{c} 0\\ \vdots \\ 0\\ tc_{g-l+1,g}^{\bar{l}}\\ 0\\ \vdots\\ 0 \end{array}\right)\\
&= \left\{\begin{array}{ll} at, & \text{if $g-l+1=f+i$,}\\ 0, & \text{otherwise}.\end{array}\right. 
\end{align*}
Similarly,
\begin{align*}
(a_{(i,\bar{j},a,s)}^{(k,\bar{l},b,t)})_f &=\left\{\begin{array}{ll}\delta_{i+f,k}ab+\delta_{j-f-1,l}st, & \text{if $f\in \{2,\ldots, n-1\}$}\\ 0, & \text{if  $f\in \{\bar{1},\ldots,\overline{n-1}\}$,}\end{array}\right.\\
(b_{(i,\bar{j},a,s)}^{(k,\bar{l},b,t)})_g&=\left\{\begin{array}{ll} \delta_{k-g,i} ab, &\text{if $g\in \{2,\ldots,n-1\}$,}\\ \delta_{g-l+1,i} at, & \text{if $g\in \{\bar{1},\ldots,\overline{n-1}\}$.}\end{array}\right. 
\end{align*}
Note that if $b_{(i,\bar{j},a,s)}^{(k,\bar{l},b,t)}$ has a nonzero entry in row $g$, then $M_{(i,\bar{j},a,s)}^{(k,\bar{l},b,t)}$ has no nonzero entry in column $g$, so in order for $b_{(i,\bar{j},a,s)}^{(k,\bar{l},b,t)}$ to be perpendicular to $\Null(M_{(i,\bar{j},a,s)}^{(k,\bar{l},b,t)})$, we have $b_{(i,\bar{j},a,s)}^{(k,\bar{l},b,t)}=0$.

Thus, $(i,\bar{j},a,s)$ meshes with $(k,\bar{l},b,t)$ if and only if 
\begin{enumerate}
\item[(1)] $k-i\geq 2$ implies $a=0$ or $b=0$,
\item[(2)] $n-1\geq i+l-1$ implies $a=0$ or $t=0$,
\item[(3)] $j-l-1\geq 2$ implies $s=0$ or $t=0$.
\end{enumerate}
In other words, $(i,\bar{j},a,s)$ meshes with $(k,\bar{l},b,t)$ if $a_{(i,\bar{j},a,s)}^{(k,\bar{l},b,t)}=b_{(i,\bar{j},a,s)}^{(k,\bar{l},b,t)}=0$ (and in this case $M_{(i,\bar{j},a,s)}^{(k,\bar{l},b,t)}=0$). Therefore, our supercharacter formula is
\begin{equation*}
\overline{\chi^{(k,\bar{l},b,t)}(x_{(i,\bar{j},a,s)})}=\left\{\begin{array}{ll} q^{\max\{k-2,n-1-l\}} \theta(ab\delta_{ik}+st\delta_{jl}), & \text{if $a_{(i,\bar{j},a,s)}^{(k,\bar{l},b,t)}=b_{(i,\bar{j},a,s)}^{(k,\bar{l},b,t)}=0$,}\\ 0, & \text{otherwise,}\end{array}\right.
\end{equation*}
where if $i=\cdot$, $j=\cdot$, $\bar{k}=\cdot$, or $\bar{l}=\cdot$, then substitute $i=2$, $j=2$, $k=n-1$, or $l=n-1$, respectively, in the formula.  Note that if we restrict this formula to the subgroup $H_{n-1}$ we obtain the supercharacter formula for $H_{n-1}$ computed in \cite{DI06}.

\subsubsection*{2. A supercharacter table}
  
The following example was used for several examples (and counterexamples) in \cite{DI06}.   Let
$$H=\left\{\left(\begin{array}{cccc} 1 & a & b & d\\ 0 & 1 & a & c\\ 0 & 0 & 1 & a\\ 0 & 0 & 0 & 1\end{array}\right)\ \bigg|\ a,b,c,d\in \FF_2\right\},$$
which is a group of order 16 with a presentation given by
$$\langle x,r,z\ \mid\ x^4=r^2=z^2=1, [x,r]=z, [x,z]=1, [r,z]=1\rangle,$$
where
$$x=\left(\begin{array}{cccc} 1 & 1 & 0 & 0\\ 0 & 1 & 1 & 0\\ 0 & 0 & 1 & 1\\ 0 & 0 & 0 & 1\end{array}\right),\quad r=\left(\begin{array}{cccc} 1 & 0 & 0 & 0\\ 0 & 1 & 0 & 1\\ 0 & 0 & 1 & 0\\ 0 & 0 & 0 & 1\end{array}\right),\quad z=\left(\begin{array}{cccc} 1 & 0 & 0 & 1\\ 0 & 1 & 0 & 0\\ 0 & 0 & 1 & 0\\ 0 & 0 & 0 & 1\end{array}\right).$$
Note that
$$l=\left(\begin{array}{cccc} 1 & 0 & 1 & 0\\ 0 & 1 & 0 & 0\\ 0 & 0 & 1 & 0\\ 0 & 0 & 0 & 1\end{array}\right)=x^2r.$$
(We use a slightly different set of generators than \cite{DI06} for $H$).  

The character table of $H$ is
$$\begin{array}{|c|cccccccccc|}\hline & \{1\} & \{z\} & \{lr\} & \{lrz\} & \{x,xz\} & \{xlr,xlrz\} & \{r,rz\} & \{l, lz\} & \{xr, xrz\} & \{xl,xlz\}\\ \hline
\chi_1 & 1 & 1 & 1 & 1 & 1 & 1 & 1 & 1 & 1 & 1\\
\chi_2 & 1 & 1 & -1 & -1 & i & -i & 1 & -1 & i & -i\\
\chi_3 & 1 & 1 & 1 & 1 & -1 & -1 & 1 & 1 & -1 & -1\\
\chi_4 & 1 & 1 & -1 & -1 & -i & i & 1 & -1 & -i & i\\
\chi_5 & 1 & 1 & 1 & 1 & 1 & 1 & -1 & -1 & -1 & -1\\
\chi_6 & 1 & 1 & -1 & -1 & i & -i & -1 & 1 & -i & i\\
\chi_7 & 1 & 1 & 1 & 1 & -1 & -1 & -1 & -1 & 1 & 1\\
\chi_8 & 1 & 1 & -1 & -1 & -i & i & -1 & 1 & i & -i\\
\chi_9 & 2 & -2 & 2 & -2 & 0 & 0 & 0 & 0 & 0 & 0\\
\chi_{10} & 2 & -2 & -2 & 2 & 0 & 0 & 0 & 0 & 0 & 0\\ \hline
\end{array}$$
The supercharacter table of $H$ is
$$\begin{array}{|c|ccccccc|}\hline & \{1\} & \{z\} & \{lr,lrz\} & \{x,xz,xlr,xlrz\} & \{r,rz\} & \{l, lz\} & \{xr, xrz,xl,xlz\}\\ \hline
\chi_1 & 1 & 1 & 1 & 1 & 1 & 1 & 1\\ 
\chi_x & 1 & 1 & 1 & -1 & 1 & 1 & -1\\
\chi_{lr} & 1 & 1  & 1 & 1 & -1 & -1 & -1\\
\chi_{rxl} & 1 & 1 & 1 & -1 & -1 & -1 & 1\\
\chi_{r} & 2 & 2 & -2 & 0 & -2 & 2 & 0\\
\chi_l & 2 & 2 & -2 & 0 & 2 & -2 & 0\\
\chi_z & 4 & -4 & 0 & 0 & 0 & 0 & 0\\ \hline
\end{array}$$
so 
$$\chi_1=\chi_1,\quad\chi_x=\chi_3,\quad\chi_{lr}=\chi_5,\quad\chi_{rxl}=\chi_7,$$
$$\chi_{r}=\chi_6+\chi_8,\quad\chi_l=\chi_2+\chi_4, \quad\chi_z=\chi_9+\chi_{10}.$$

There are many more supercharacter theories (in the sense of \cite{DI06}) and we can find many by inspection, such as
$$\begin{array}{|c|cccccc|}\hline & \{1\} & \{z\} & \{lr,lrz\} & \{x,xz,xlr,xlrz, xr, xrz,xl,xlz\} & \{r,rz\} & \{l, lz\}\\ \hline
\chi_1 & 1 & 1 & 1 & 1 & 1 & 1 \\ 
\chi_x & 1 & 1 & 1 & -1 & 1 & 1 \\
\tilde\chi_{lr} & 2 & 2  & 2 & 0 & -2 & -2 \\
\chi_{r} & 2 & 2 & -2 & 0 & -2 & 2 \\
\chi_l & 2 & 2 & -2 & 0 & 2 & -2 \\
\chi_z & 4 & -4 & 0 & 0 & 0 & 0 \\ \hline
\end{array}$$
with $\tilde\chi_{lr}=\chi_{lr}+\chi_{rxl}$;
$$\begin{array}{|c|cccc|}\hline & \{1\} & \{z\} & \{l,r,lr,lz,rz,lrz\} & \{x,xz,xlr,xlrz, xr, xrz,xl,xlz\}\\ \hline
\chi_1 & 1 & 1 & 1 & 1  \\ 
\chi_x & 1 & 1 & 1 & -1  \\
\tilde\chi_l & 6 & 6 & -2 & 0  \\
\chi_z & 4 & -4 & 0 & 0 \\ \hline
\end{array}$$
with $\tilde\chi_l=\chi_l+\chi_r+\tilde\chi_{lr}$;  or
$$\begin{array}{|c|ccc|}\hline & \{1\} & \{z\} & \{l,r,lr,lz,rz,lrz, x,xz,xlr,xlrz, xr, xrz,xl,xlz\}\\ \hline
\tilde\chi_1 & 1 & 1 & 1   \\ 
\tilde\chi_x & 7 & 7 & -1   \\
\tilde\chi_z & 4 & -4 & 0  \\ \hline
\end{array}$$
with $\tilde{\chi}_x=\tilde\chi_l+\chi_x$.

\end{document}